\documentclass[reqno]{amsart}
\usepackage{mathrsfs}
\usepackage{color}
\usepackage{amsmath}
\usepackage{amsfonts}
\usepackage{amssymb}
\usepackage{graphicx}%
\usepackage{hyperref}

 \newtheorem{Theorem}{Theorem}[section]
 \newtheorem{Corollary}[Theorem]{Corollary}
 \newtheorem{Lemma}[Theorem]{Lemma}
 \newtheorem{Proposition}[Theorem]{Proposition}

\newtheorem{Question}[Theorem]{Question}

 \newtheorem{Remark}[Theorem]{Remark}

 \numberwithin{equation}{section}

\begin{document}

\title[effectiveness of strong openness conjecture]
 {effectiveness of Demailly's strong openness conjecture and related problems}

\author{Qi'an Guan}
\address{Qi'an Guan:  Beijing International Center for Mathematical Research, and School of Mathematical Sciences,
Peking University, Beijing, 100871, China.}
\email{guanqian@math.pku.edu.cn}
\author{Xiangyu Zhou}
\address{Xiangyu Zhou: Institute of Mathematics, AMSS, and Hua Loo-Keng Key Laboratory of Mathematics, Chinese Academy of Sciences, Beijing, China}
\email{xyzhou@math.ac.cn}

\thanks{The authors were partially supported by NSFC}

\subjclass[2010]{32D15, 32E10, 32L10, 32U05, 32W05}

\keywords{strong openness conjecture, effectiveness,
plurisubharmonic function, multiplier ideal sheaf, complex
singularity exponent}

\date{}

\dedicatory{}

\commby{}


\begin{abstract}
In this article, stimulated by the effectiveness in Berndtsson's solution of the openness conjecture
and continuing our solution of Demailly's strong openness conjecture,
we discuss conditions to guarantee the effectiveness of the conjecture
and establish such an effectiveness result.
We explicitly point out a lower semicontinuity property of plurisubharmonic functions with a multiplier,
which is implicitly contained in \cite{GZopen-b}.
We also obtain optimal effectiveness of the conjectures of
Demailly-Koll\'{a}r and Jonsson-Mustat\u{a} respectively.
\end{abstract}

\maketitle

\section{Introduction}

\subsection{Background: Strong openness conjecture}\label{sub:strong}
$\\$

In
\cite{GZopen-a,GZopen-b}, the authors solved
the strong openness conjecture posed by Demailly in \cite{demailly-note2000} and \cite{demailly2010}
(see also \cite{DEL00}, \cite{D-P03}, \cite{BFM08},
\cite{Leh11}, \cite{JM12}, \cite{Gue12}, \cite{Cao12},
\cite{JM13}, \cite{Yi2013}, \cite{Mat2013}, etc. ):

\emph{\textbf{Strong openness conjecture:}  Let $\varphi$ be a plurisubharmonic function on a complex manifold.
Then $$\mathcal{I}_{+}(\varphi)=\mathcal{I}(\varphi),$$
where $\mathcal{I}(\varphi)$ is the multiplier ideal sheaf and $\mathcal{I}_{+}(\varphi):=\cup_{\varepsilon>0}\mathcal{I}((1+\varepsilon)\varphi)$.}

Recall that the multiplier ideal sheaf $\mathcal{I}(\varphi)$ is the sheaf of germs of holomorphic functions $f$ such that
$|f|^{2}e^{-\varphi}$ is locally integrable (see \cite{Nadel90}, see also \cite{siu05},
\cite{siu09}, \cite{demailly2010}, etc.).

For $dim X \leq 2$, the strong openness conjecture was
proved in \cite{JM12} by studying the asymptotic jumping numbers for
graded sequences of ideals.

It is not hard to see that the truth of
the strong openness conjecture is equivalent to the following theorem:

\begin{Theorem}
\label{t:strong1015}\cite{GZopen-a,GZopen-b}
Let $\varphi$ be a negative plurisubharmonic function on the unit polydisc $\Delta^{n}\subset\mathbb{C}^{n}$,
which satisfies
$$\int_{\Delta^{n}}|F|^{2}e^{-\varphi}d\lambda_{n}<+\infty,$$
where  $d\lambda_{n}$
is the Lebesgue measure on $\mathbb{C}^{n}$,
$F$ is a holomorphic function on $\Delta^{n}$. Then there
exists a number $p>1$, such that
$$\int_{\Delta^{n}_{r}}|F|^{2}e^{-p\varphi}d\lambda_{n}<+\infty,$$
where $r\in(0,1)$.
\end{Theorem}

Assuming $\mathcal{I}(\varphi)=\mathcal{O}_{X}$, the strong openness conjecture degenerates to the openness conjecture posed by Demailly and Koll\'{a}r in \cite{D-K01}.
The dimension two case of the Openness conjecture was proved by Favre and Jonsson in \cite{FM05j} (see also \cite{FM05v}).

Recently, Berndtsson \cite{berndtsson13} proved the openness conjecture.
Actually, Berndtsson also obtained effectiveness in his solution of the openness conjecture.

Therefore it is natural to ask a similar effectiveness problem in the case of the strong openness conjecture:

\begin{Question}
\label{q:effect}
Under what kinds of conditions,
one can find effective $p>1$, such that
$$(F,z_{0})\in\mathcal{I}(p\varphi)_{z_{0}}.$$
\end{Question}

In the following subsection, we will discuss Question \ref{q:effect} and find suitable conditions.

\subsection{Effectiveness of the strong openness conjecture}
$\\$

Let $D$ be a given pseudoconvex domain in $\mathbb{C}^{n}$, and $z_{0}$ be a point in $D$.
Let $\varphi$ be a negative plurisubharmonic function on $D$,
and $F$ be a holomorphic function on $D$.

Denote by
$$\parallel F\parallel_{\varphi}:=(\int_{D}|F|^{2}e^{-\varphi}d\lambda_{n})^{1/2}.$$

Let $C_{1}$ be a positive constant.

First natural condition is:

1) $\parallel F\parallel_{\varphi}^{2}\leq C_{1}.$

The following example shows that under this condition there doesn't
exist an effective (only depending on $C_{1}$) number $p>1$,
such that $(F,z_{0})\in\mathcal{I}(p\varphi)_{z_{0}}$:

For $n=1$ case, let $D=\Delta$, $z_{0}=0$, and $F=z^{m}$, and $\varphi=2m\log|z|$. Then $p$ must be smaller than $1+\frac{1}{m}$.
Note that $\int_{\Delta}|F|^{2}e^{-\varphi}d\lambda_{n}=\pi$ for any $m$. Thus $p$ depends on $F$.

Thus the condition $1)$ is not enough to answer Question \ref{q:effect}.
$\\$

Then we consider the following modified conditions:

1) $\parallel F\parallel_{\varphi}^{2}\leq C_{1}$;

2) given $F$.

The following example shows that under this condition there still doesn't
exist an effective number $p>1$ (only depending on $C_{1}$ and $F$),
such that $(F,z_{0})\in\mathcal{I}(p\varphi)_{z_{0}}$:

For $n=2$ case, let $D=\mathbb{B}^{2}$,
$z_{0}=(0,0)$,
$F=z_{1}$,
and $\varphi_{\theta,\delta}=2\delta\log(|z_{1}\cos\theta+z_{2}\sin\theta|)$.
Note that for any $\delta\in(0,1)$, there exists $\theta=\theta(\delta)>0$ small enough,
such that
$$\int_{\Delta^{2}}|F|^{2}e^{-\varphi_{\theta,\delta}}d\lambda_{n}<2\pi^{2}.$$
However,
$$\int_{\Delta^{2}_{r}}|F|^{2}e^{-\frac{1}{\delta}\varphi_{\theta,\delta}}d\lambda_{n}=+\infty$$
for any $r\in(0,1)$.
When $\delta$ go to zero,
$p$ must goes to 1.

Thus the modified conditions $1)$ and $2)$ are not enough to answer Question \ref{q:effect}.
$\\$

Therefore,
in order to answer Question \ref{q:effect},
we need to consider the pair $(F,\varphi)$ instead of $F$ and $\varphi$ separately.

We generalize the Bergman kernel $K$ on $D$ to $K_{\varphi,F}$ as follows:

$$K_{\varphi,F}(z_{0}):=\frac{1}{\inf\{\|F_{1}\|^{2}_{0}|(F_{1}-F,z_{0})\in\mathcal{I}_{+}(2c_{z_{0}}^{F}(\varphi)\varphi)_{z_0}{\,}\&{\,}F_{1}\in\mathcal{O}(D)\}}$$
where $c_{z_{0}}^{F}(\varphi):=\sup\{c\geq0:|F|^{2}e^{-2c\varphi}$ is $L^1$ on a neighborhood of $z_{0}\}$ is the jumping number (see \cite{JM13}).
Especially, when $F\equiv 1$, $c_{z_{0}}^{F}(\varphi)$ will degenerate to
the complex singularity exponent $c_{z_{0}}(\varphi)$ (or log canonical threshold) (see \cite{tian87,Sho92,Ko92,D-K01}, etc.).

It is clear that when $F\equiv1$ and $\varphi$ is the pluricomplex Green function $G(z,z_{0})$ on $D$
(i.e the upper envelop of negative plurisubharmonic functions on $D$,
which satisfy $G_{z,w}-n\log|z-w|^{2}$ is locally finite near $w$, see \cite{demailly86}),
then
$K_{\varphi,F}(z_{0})$ degenerates to the Bergman kernel $K(z_{0}).$
$\\$

Let $D$ be $\Delta\subset\mathbb{C}$.
When $F=z^{m}$ and $\varphi=\log|z|^{2m}$,
then
$$K_{\varphi,F}(o)=(\int_{\Delta}|z|^{2m}\lambda_{1})^{-1}=\frac{m+1}{\pi}.$$

Let $D$ be $\Delta^{2}\subset\mathbb{C}^{2}$.
When $F=z_{1}\cos\theta+z_{2}\sin\theta$ and $\varphi=2\delta\log|z_{1}|$,
then
$$K_{\varphi,F}(o)=(\sin^{2}\theta\int_{\Delta^{2}}|z_{2}|^{2}\lambda_{2})^{-1}=\frac{2}{\pi^{2}\sin^{2}\theta}.$$

Let $D$ be $\Delta^{2}\subset\mathbb{C}^{2}$.
When $F=z_{1}+z^{2}_{2}$ and $\varphi=2\log|z_{1}|$,
then
$$K_{\varphi,F}(o)=(\int_{\Delta^{2}}|z_{2}|^{4}\lambda_{2})^{-1}=\frac{4}{\pi^{2}}.$$
$\\$

We define a useful function to establish the effectiveness of the strong openness conjecture:
$$\theta(t):=(\frac{1}{(t-1)(2t-1)})^{\frac{1}{t}},$$
where $t\in(1,+\infty)$.
$\\$

In the present article,
we establish the effectiveness of the strong openness conjecture as follows:

\begin{Theorem}
\label{t:effect}
Let $C_{1}$ and $C_{2}$ be two positive constants.
We consider the set of the pairs $(F,\varphi)$ satisfying

$1)$ $\parallel F\parallel^{2}_{\varphi}\leq C_{1}$;

$2)$ $K^{-1}_{\varphi,F}(z_{0})\geq C_{2}$.

Then for any $p>1$ satisfying
$$\theta(p)>\frac{C_{1}}{C_{2}},$$
we have
$$(F,z_{0})\in\mathcal{I}(p\varphi)_{z_{0}}.$$
\end{Theorem}

Note that the proof of Theorem \ref{t:effect} does not depend on the truth of the strong openness conjecture,
then
Theorem \ref{t:effect} can be regarded as a presentation of our proof of the strong openness conjecture.

Especially, letting $F\equiv 1$ in Theorem \ref{t:effect},
noting that
$$K_{\varphi,1}(z_{0})\leq K(z_{0}),$$
one can obtain Berndtsson's effectiveness of the openness conjecture in \cite{berndtsson13}:
\begin{Corollary}
\label{c:effect}
Let $C_{1}$ and $C_{2}$ be two positive constants.
We consider the set of $\varphi$ satisfying

$1)$ $\parallel 1\parallel^{2}_{\varphi}=\int_{D}e^{-\varphi}d\lambda_{n}\leq C_{1}$;

$2)$ $K^{-1}(z_{0})\geq C_{2}$.
$\\$
Then for any $p>1$ satisfying
\begin{equation}
\label{equ:sharp20140128d}
\theta(p)> \frac{C_{1}}{C_{2}},
\end{equation}
we have $e^{-p\varphi}$ is integrable near $z_{0}$.
\end{Corollary}

When $D=\mathbb{B}^{n}$,
recall that $\varepsilon_{0}$ in \cite{berndtsson13},
$$\varepsilon_{0}:=\inf\{\|w\|^{2}_{0}|\exists z_{0}\in\mathbb{B}^{n}_{1/2}, |w(z_{0})|=\frac{1}{10}\}=\inf_{z_{0}\in\mathbb{B}^{n}_{1/2}}\inf\{\|w\|^{2}_{0}||w(z_{0})|=\frac{1}{10}\},$$
and
$$K(z_{0})=\sup_{w\in\mathcal{O}_{\mathbb{B}^{n}}}\frac{|w(z_{0})|^{2}}{\|w\|^{2}_{0}}
=\sup_{|w(z_{0})|=\frac{1}{10}}\frac{|w(z_{0})|^{2}}{\|w\|^{2}_{0}}=\frac{\frac{1}{100}}{\inf_{|w(z_{0})|=\frac{1}{10}}\|w\|^{2}_{0}}.$$
It follows that
\begin{equation}
\label{equ:bernd}
\varepsilon_{0}=\inf_{z_{0}\in\mathbb{B}^{n}_{1/2}}\frac{1}{100}\frac{1}{K(z_{0})}.
\end{equation}

The following Remark tells us that Corollary \ref{c:effect}
is (a more precise version of) Berndtsson's effectiveness result of the openness conjecture:
\begin{Remark}
\label{r:effect}
When $D=\mathbb{B}^{n}$,
and $|z_{0}|\leq\frac{1}{2}$,
Berndtsson in \cite{berndtsson13} pointed out that
if $p$
satisfies
$$\frac{1}{(p-1)}\geq\parallel 1\parallel^{2}_{\varphi}\frac{2}{\varepsilon_{0}},$$
then $e^{-p\varphi}$ is integrable near $z_{0}$.

One can check that
\begin{equation}
\label{equ:sharp20140126}
\frac{1}{200(t-1)}<\frac{1}{6(t-1)}<\frac{t}{6(t-1)}<(\frac{1}{(t-1)(2t-1)})^{\frac{1}{t}}=\theta(t)
\end{equation}
for any $t\in(1,+\infty)$.

By inequality \ref{equ:sharp20140128d}, \ref{equ:bernd} and \ref{equ:sharp20140126}, it follows that
Corollary \ref{c:effect}
is (a more precise version of) Berndtsson's effectiveness of the openness conjecture.
\end{Remark}

In subsection \ref{sec:2014springb},
we show that
\begin{equation}
\label{equ:sharp20140127a}
\frac{1}{\sqrt{3}e^{\frac{2}{e}}}\frac{t}{(t-1)}<(\frac{1}{(t-1)(2t-1)})^{\frac{1}{t}}.
\end{equation}

In subsection \ref{sec:2014spring},
we give a more precise form of inequality \ref{equ:sharp20140127a}.
$\\$

We establish Theorem \ref{t:effect} by the following Proposition:

\begin{Proposition}
\label{p:effect}
Let $C_{1}$ and $C_{2}$ be two positive constants.
We consider the set of the triples $(F,\varphi,p)$
satisfying

$1)$ $\parallel F\parallel^{2}_{D,\varphi}\leq C_{1}$;

$2)$ $C_{F,p\varphi}(z_{0})\geq C_{2}$,
$\\$
where $p>1$ is a real number
and $C_{F,p\varphi}(z_{0})$ is the infimum of $\int_{D}|F_{1}|^{2}d\lambda_{n}$ for all $F_{1}\in\mathcal{O}(D)$
satisfying condition $(F_{1}-F,z_{0})\in\mathcal{I}(p\varphi)_{z_{0}}$.

Then $p$ must satisfy
\begin{equation}
\label{equ:sharp20140128a}
\theta(p)\leq \frac{C_{1}}{C_{2}}.
\end{equation}
\end{Proposition}

\subsection{A lower semicontinuity property of plurisubharmonic functions with a multiplier}\label{sub:lower}
$\\$

In \cite{D-K01},
Demailly and Koll\'{a}r conjectured that:

For every nonzero holomorphic function $f$ on $X$,
there is a number $\delta=\delta(f,K,L)>0$, such that for any holomorphic function
$g$ on $X$ with
$$\sup_{L}|g-f|<\delta\Rightarrow c_{K}(\log|g|)\geq c_{K}(\log|f|),$$
where the compact set $K$ contained in an open subset $L$ of complex manifold $X$.

In \cite{D-K01},
the authors proved that
the above conjecture
is implied by the ACC conjecture (see \cite{Sho92} or \cite{Ko92}).
The ACC conjecture was proved by Hacon, McKernan and Xu in \cite{HKX2012}.

Note that $c\log|f|$ is a plurisubharmonic function,
then
we \cite{GZopen-b} replaced $c\log|f|$ by general plurisubharmonic functions,
and
obtained the following lower semicontinuity property of plurisubharmonic functions,
which is a new proof of the above conjecture in \cite{D-K01} without using the ACC conjecture:

\begin{Proposition}
\label{p:lower_open}\cite{GZopen-b}
Let $\{\phi_{m}\}_{m=1,2,\cdots}$ be a sequence of negative plurisubharmonic functions on $\Delta^{n}$,
which is convergent to a negative Lebesgue measurable function $\phi$ on $\Delta^{n}$ in Lebesgue measure.
Assume that $e^{-\phi_{m}}$ are all not integrable near $o$.
Then $e^{-\phi}$ is not integrable near $o$.
\end{Proposition}

In fact, our proof of Proposition \ref{p:lower_open} in \cite{GZopen-b} already contains the following lower semicontinuity property of
plurisubharmonic functions on multiplier ideal sheaves:

\begin{Proposition}
\label{p:lower}
Let $\{\phi_{m}\}_{m=1,2,\cdots}$ be a sequence of negative plurisubharmonic functions on $\Delta^{n}$,
which is convergent to a negative Lebesgue measurable function $\phi$ on $\Delta^{n}$ in Lebesgue measure.
Let $\{F_{m}\}_{m=1,2,\cdots}$ be a sequence of holomorphic functions on $\Delta^{n}$ with uniform bound,
which is convergent to a Lebesgue measurable function $F$ on $\Delta^{n}$ in Lebesgue measure.
Assume that for any neighborhood $U$ of $o$,
the pairs $(F_{m},\phi_{m})$ $(m=1,2,\cdots)$ satisfying
$$\inf_{m}K^{-1}_{\phi_{m},F_{m}}(o)>0,$$
where $K_{\phi{m},F_{m}}$ is the generalized Bergman kernel on $U$.
Then
$|F|^{2}e^{-\varphi}$ is not integrable near $o$.
If $\phi$ is plurisubharmonic
and $F$ is holomorphic,
then
$$(F,o)\not\in\mathcal{I}(\phi)_{o}.$$
\end{Proposition}

Especially, letting $F\equiv 1$ in Proposition \ref{p:lower},
we obtain Proposition \ref{p:lower_open}.

In fact, Proposition \ref{p:lower} can be regarded as another presentation of our proof of the strong openness conjecture.

An another presentation of Proposition \ref{p:lower} is:

Let $\{\phi_{m}\}_{m=1,2,\cdots}$ be a sequence of negative plurisubharmonic functions on $\Delta^{n}$,
which is convergent to a plurisubharmonic function $\phi$ on $\Delta^{n}$ in Lebesgue measure.

Let $\{F_{m}\}_{m=1,2,\cdots}$ be a sequence of holomorphic functions on $\Delta^{n}$ with uniform bound,
which is convergent to a holomorphic function $F$ on $\Delta^{n}$ in Lebesgue measure.

Assume that for any neighborhood $U$ of $o$,
the pairs $(F_{m},\phi_{m})$ $(m=1,2,\cdots)$ satisfying
$$\inf_{m}K^{-1}_{\phi_{m},F_{m}}(o)>0,$$

There exists $m_{0}$, such that for any $m\geq m_{0}$,
$$c_{o}^{F_{m}}(\varphi_{m})\geq c_{o}^{F}(\varphi),$$
where
$c_{o}^{F}(\psi)=sup\{c\geq0:|F|^{2}e^{-2c\psi}$ is $L^1$ on a neighborhood of $o\}$ is the jumping number (see \cite{JM13}).

\subsection{Optimal effectiveness of a conjecture posed by Demailly and Koll\'{a}r}
$\\$
In \cite{GZopen-b}, we solved the following conjecture
about the volume growth of the
sublevel sets of plurisubharmonic functions related to the complex singularity
exponents posed by Demailly and Koll\'{a}r
in \cite{D-K01} (see also \cite{FM05j}, \cite{FM05v}, \cite{JM12} and \cite{JM13}, etc.):

\emph{\textbf{Conjecture D-K:} Let $\varphi$ be a plurisubharmonic function on $\Delta^{n}\subset\mathbb{C}^{n}$,
and $K$ be compact subset of $\Delta^{n}$.
If $c_{K}(\varphi)<+\infty$,
then
$$\frac{1}{r^{2c_{K}(\varphi)}}\mu(\{\varphi<\log r\})$$
has
a uniform positive lower bound independent of $r\in(0,1)$ small enough,
where $c_{K}(\varphi)=sup\{c\geq0:e^{-2c\varphi}$ is $L^1$ on a neighborhood of $K\}$,
and $\mu$ is the Lebesgue volumes on $\mathbb{C}^{n}$.}

For $n\leq 2$, the above conjecture
was proved by Favre and Jonsson in \cite{FM05j} (see also \cite{FM05v}).

In \cite{GZopen-b}, in order to prove Conjecture D-K, we obtained an estimate about the volume growth of the sublevel sets
of plurisubharmonic functions:

\begin{Theorem}
\label{t:GZ_JM}\cite{GZopen-b}
Let $\varphi$ be a plurisubharmonic function on $\Delta^{n}\subset\mathbb{C}^{n}$.
Let $F$ be a holomorphic function on $\Delta^{n}$.
Assume that $|F|^{2}e^{-\varphi}$ is not locally integrable near $o$.
Then
$$\int_{\Delta^{n}}\mathbb{I}_{\{-(R+1)<\varphi<-R\}}|F|^{2}e^{-\varphi}d\lambda_{n}$$
has a uniform positive lower bound independent of $R>>0$.
Especially, if $F=1$, then
$$e^{R}\mu(\{-(R+1)<\varphi<-R\})$$
has a uniform positive lower bound independent of $R>>0$.
\end{Theorem}

Theorem \ref{t:GZ_JM} tells us that
$$\liminf_{R\to+\infty}\int_{\Delta^{n}}\mathbb{I}_{\{-(R+1)<\varphi<-R\}}|F|^{2}e^{-\varphi}d\lambda_{n}>0.$$

In fact, our proof of Theorem \ref{t:GZ_JM} in \cite{GZopen-b} already contains the following
effectiveness of the uniform positive lower bound in Theorem \ref{t:effect}:

\begin{Proposition}
\label{p:effect_conj.D-K}Let $B_{0}\in(0,1]$ be arbitrarily given.
Let $\varphi$ be a negative plurisubharmonic function on pseudoconvex domain $D\subset\mathbb{C}^{n}$.
Let $F$ be a holomorphic function on $D$.
Assume that $|F|^{2}e^{-\varphi}$ is not locally integrable near $z_{0}$.
Then we obtain that
$$\liminf_{R\to+\infty}\frac{1}{B_{0}}\int_{D}\mathbb{I}_{\{-(R+B_{0})<\varphi<-R\}}|F|^{2}e^{t_{0}+B_{0}}d\lambda_{n}\geq K_{\varphi,F}^{-1}(z_{0})$$
Especially, if $F=1$, then
\begin{equation}
\label{equ:DK.GZ.optimal}
\liminf_{R\to+\infty}e^{R+B_{0}}\frac{1}{B_{0}}\mu(\{-(R+B_{0})<\varphi<-R\})\geq K_{\varphi,1}^{-1}(z_{0})\geq K^{-1}(z_{0}).
\end{equation}
\end{Proposition}

Taking $R=kB_{0}$ in inequality \ref{equ:DK.GZ.optimal},
for any given $\varepsilon>0$, there exists $k_{0}$ depending on $B_{0}$, such that for any $k\geq k_{0}$,
on can obtain
$$e^{(k+1)B_{0}}\frac{1}{B_{0}}\mu(\{-(k+1)B_{0}<\varphi<-kB_{0}\})\geq (K^{-1}(z_{0})-\varepsilon),$$
i.e.,
\begin{equation}
\label{equ:DK.20140301a}
\mu(\{-(k+1)B_{0}<\varphi<-kB_{0}\})\geq e^{-(k+1)B_{0}}B_{0}(K^{-1}(z_{0})-\varepsilon).
\end{equation}

Taking sum with $k\geq k_{0}$ in inequality \ref{equ:DK.20140301a} and letting $B_{0}$ goes to $0$,
one can obtain
$$\liminf_{R\to +\infty}\mu(\{\varphi<-R\})e^{R}\geq (K^{-1}(z_{0})-\varepsilon).$$

Then we obtain an optimal estimate of the lower bound of  $\liminf_{R\to+\infty}e^{R}\mu(\{\varphi<-R\})$
\begin{equation}
\label{equ:DK.optimal}
\liminf_{R\to+\infty}e^{R}\mu(\{\varphi<-R\})\geq K_{\varphi,1}^{-1}(z_{0})\geq K^{-1}(z_{0}).
\end{equation}

When $D=\Delta\subset\mathbb{C}$, $\varphi=\log|z|^{2}$ and $z_{0}=0$, the equality in inequality \ref{equ:DK.GZ.optimal} holds.

Replacing $R$ by $-2c_{K}(\varphi)\log r$ and $\varphi$ by $2c_{K}(\varphi)\varphi$,
we obtain the optimal effectiveness of Conjecture D-K:
$$\liminf_{r\to 0}\frac{1}{r^{2c_{K}(\varphi)}}\mu(\{\varphi<\log r\})\geq K^{-1}(z_{0})\geq \inf_{z\in K}K^{-1}(z),$$
where $z_{0}$ satisfies $c_{z_{0}}(\varphi)=c_{K}(\varphi)$.

When $D=\Delta\subset\mathbb{C}$, $\varphi=\log|z|^{2}$ and $K=\{0\}$, the equality in the above inequality holds.

\subsection{Optimal effectiveness of a conjecture posed by Jonsson and Mustat\u{a}}
$\\$
In \cite{GZopen-b},
we solved the following conjecture about the volumes growth of the
sublevel sets of quasi-plurisubharmonic functions posed by Jonsson and Mustat\u{a}
in \cite{JM13} (see also \cite{JM12}):

\emph{\textbf{Conjecture J-M:} Let $\psi$ be a plurisubharmonic function on $\Delta^{n}\subset\mathbb{C}^{n}$.
If $c_{o}^{I}(\psi)<+\infty$,
then
$$\frac{1}{r^2}\mu(\{c^{I}_{o}(\psi)\psi-\log|I|<\log r\})$$
has
a uniform positive lower bound independent of $r\in(0,1)$ small enough,
where
$I$ is an ideal of $\mathcal{O}_{\Delta^{n},o}$, which is generated by $\{f_{j}\}_{j=1,\cdots,l}$,
$$\log|I|:=\log\max_{1\leq j\leq l}|f_{j}|,$$
$c_{o}^{I}(\psi)=sup\{c\geq0:|I|^{2}e^{-2c\psi}$ is $L^1$ on a neighborhood of $o\}$ is the jumping number in \cite{JM13}.}

For $n\leq 2$, the above conjecture
was proved by Jonsson and Mustat\u{a} in \cite{JM12}.

In \cite{GZopen-b}, in order to prove conjecture J-M, we gave the following estimate about the volume growth of the sublevel sets
of quasiplurisubharmonic functions:

\begin{Theorem}
\label{t:GZ_JM201312}\cite{GZopen-b}
Let $\psi$ be a plurisubharmonic function on $\Delta^{n}$,
and $F$ be a holomorphic function on $\Delta^{n}$.
Assume that $|F|^{2}e^{-\psi}$ is not locally integrable near $o$.
Then
$$e^{R}\frac{1}{B_{0}}\mu(\{-R-B_{0}<\psi-\log|F|^{2}<-R\})$$
has
a uniformly positive lower bound independent of $R>>0$ and $B_{0}\in(0,1]$.
\end{Theorem}

Theorem \ref{t:GZ_JM201312} tells us that
$$\liminf_{R\to+\infty}e^{R}\frac{1}{B_{0}}\mu(\{-R-B_{0}<\psi-\log|F|^{2}<-R\})>0$$

In fact, our proof of Theorem \ref{t:GZ_JM201312} in \cite{GZopen-b} already contains the following
effectiveness of the uniform positive lower bound in Theorem \ref{t:effect}:

\begin{Proposition}
\label{p:effect_conj.J-M}Let $\delta$ be an arbitrarily given positive integer.
Let $\psi$ be a bounded from above plurisubharmonic function on pseudoconvex domain $D\subset\mathbb{C}^{n}$.
Let $F$ be a bounded holomorphic function on $D$.
Assume that $|F|^{2}e^{-\psi}$ is not locally integrable near $z_{0}$.
Then we obtain that
\begin{equation}
\label{equ:JM20140301}
\begin{split}
&\liminf_{R\to+\infty}e^{R}\frac{1}{B_{0}}\mu(\{-R-B_{0}<\psi-\log|F|^{2}<-R\})
\geq \frac{C_{\psi,F,\delta}}{(1+\frac{1}{\delta})e^{B_{0}}}.
\end{split}
\end{equation}
where $C_{\psi,F,\delta}=\frac{K_{\psi+\delta\max\{\psi,\log|F|^{2}\},F^{1+\delta}}^{-1}(z_{0})}{\sup_{D}e^{(1+\delta)\max\{\psi,2\log|F|\}}}$,
and $K_{\psi+\delta\max\{\psi,\log|F|^{2}\},F^{1+\delta}}$ is the generalized Bergman kernel on $D$.
\end{Proposition}

Taking $R=kB_{0}$ in inequality \ref{equ:JM20140301},
for any given $\varepsilon>0$, there exists $k_{0}$ depending on $B_{0}$, such that for any $k\geq k_{0}$,
one can obtain
$$e^{(k+1)B_{0}}\frac{1}{B_{0}}\mu(\{-(k+1)B_{0}<\psi-\log|F|^{2}<-kB_{0}\})\geq (\frac{C_{\psi,F,\delta}}{(1+\frac{1}{\delta})e^{B_{0}}}-\varepsilon),$$
i.e.
\begin{equation}
\label{equ:JM.20140301a}
\mu(\{-(k+1)B_{0}<\psi-\log|F|^{2}<-kB_{0}\})\geq e^{-(k+1)B_{0}}B_{0} (\frac{C_{\psi,F,\delta}}{(1+\frac{1}{\delta})e^{B_{0}}}-\varepsilon).
\end{equation}

Taking sum $k\geq k_{0}$ in inequality \ref{equ:JM.20140301a},
and letting $B_{0}$ go to zero,
one can obtain the following estimate:
\begin{equation}
\label{equ:2014113b}
\begin{split}
\liminf_{R\to+\infty}e^{R}\mu(\{\psi-\log|F|^{2}<-R\})
\geq \sup_{\delta\in\{1,2,\cdots\}}\frac{C_{\psi,F,\delta}}{(1+\frac{1}{\delta})}.
\end{split}
\end{equation}

By Theorem \ref{t:strong1015},
it follows that $|I|^{2}e^{-2c_{o}^{I}(\psi)\psi}$ is not integrable on any neighborhood of $o$.

Replacing $\psi$ by  $2c_{o}^{I}(\psi)\psi$, and  $R$ by $-2\log r$ in equality \ref{equ:2014113b},
we obtain the optimal effectiveness of conjecture J-M:

\begin{equation}
\label{equ:20140304}
\begin{split}
&\liminf_{r\to0}\frac{1}{r^2}\mu(\{c^{I}_{o}(\psi)\psi-\log|I|<\log r\})
\\&\geq \sup_{\delta\in\{1,2,\cdots\}}\frac{\max_{1\leq i\leq l}\{K_{2c_{o}^{I}(\psi)\psi+\delta\max\{2c_{o}^{I}(\psi)\psi,\log|f_{i}|^{2}\},f_{i}^{1+\delta}}^{-1}(o)\}}
{(1+\frac{1}{\delta})\sup_{D}e^{(1+\delta)\max\{2c_{o}^{I}(\psi)\psi,2\log|I|\}}},
\end{split}
\end{equation}
where $D$ is a relatively compact pseudoconvex domain in $\Delta^{n}$,
and
$$K_{2c_{o}^{I}(\psi)\psi+\delta\max\{2c_{o}^{I}(\psi)\psi,\log|f_{i}|^{2}\},f_{i}^{(1+\delta)}}$$ is the generalized Bergman Kernel on $D$.

\begin{Remark}
\label{r:JM-1}
When $D=\Delta$, $F=1$, $\psi=\log|z|^{2}$,
and $\delta$ goes to $\infty$, the equality in inequality \ref{equ:20140304} holds.
\end{Remark}

\begin{Remark}\label{r:JM-DK}
The optimal effectiveness of conjecture J-M is a generalization of the optimal effectiveness of conjecture D-K.
\end{Remark}

\section{Proof of effectiveness of strong openness conjecture}

\subsection{A Lemma used to prove Proposition \ref{p:effect}}
$\\$

We prove Proposition \ref{p:effect} by the following Lemma,
whose various forms already appear in \cite{guan-zhou13p,guan-zhou13ap,GZopen-b} etc. and whose proof will appear in the section 4 for the sake of completeness:

\begin{Lemma} \label{p:GZ_JM_sharp}
Let $B_{0}\in(0,1]$ be arbitrarily given.
Let $D_{v}$ be a strongly pseudoconvex domain relatively compact in
$\Delta^{n}$ containing $o$. Let $F$ be a holomorphic function on
$\Delta^{n}$. Let $\psi$ be a negative plurisubharmonic function
on $\Delta^{n}$, such that $\psi(o)=-\infty$. Then there exists a
holomorphic function $F_{v,t_{0}}$ on $D_{v}$, such that,
$$(F_{v,t_{0}}-F,o)\in\mathcal{I}(\psi)_{o}$$
and
\begin{equation}
\label{equ:3.4}
\begin{split}
&\int_{ D_v}|F_{v,t_0}-(1-b_{t_0}(\psi))F|^{2}d\lambda_{n}
\\\leq&(1-e^{-(t_{0}+B_{0})})\int_{D_v}\frac{1}{B_{0}}(\mathbb{I}_{\{-t_{0}-B_{0}<t<-t_{0}\}}\circ\psi)|F|^{2}e^{-\psi}d\lambda_{n},
\end{split}
\end{equation}
where
$b_{t_{0}}(t)=\int_{-\infty}^{t}\frac{1}{B_{0}}\mathbb{I}_{\{-t_{0}-B_{0}< s<-t_{0}\}}ds$,
and $t_{0}$ is a positive number.
\end{Lemma}

\begin{Remark}
\label{r:GZ_JM_sharp}
Replacing the strong pseudoconvexity of $D_{v}$ by pseudoconvexity, and $o$ by $z_{0}\in D_{v}$ for any $v$,
Lemma \ref{p:GZ_JM_sharp} also holds.
\end{Remark}

\subsection{Proof of Proposition \ref{p:effect}}
$\\$

As $C_{F,p\varphi}(z_{0})\geq C_{2}>0$,
then we obtain $\varphi(z_{0})=-\infty$.

Let $\psi:=p\varphi$.

Using Lemma \ref{p:GZ_JM_sharp} and Remark \ref{r:GZ_JM_sharp},
we obtain that
\begin{equation}
\label{equ:effect20140120.1}
\begin{split}
&((\int_{D_{v}}|F_{v,t_0}|^{2}d\lambda_{n})^{1/2}-(\int_{D_{v}}|(1-b_{t_0}(\psi))F|^{2}d\lambda_{n})^{1/2})^{2}
\\&\leq \int_{D_{v}}|F_{v,t_0}-(1-b_{t_0}(\psi))F|^{2}d\lambda_{n}
\\&\leq(1-e^{-(t_{0}+B_{0})})\int_{D_{v}}\frac{1}{B_{0}}(\mathbb{I}_{\{-t_{0}-B_{0}<t<-t_{0}\}}\circ\psi)|F|^{2}e^{-\psi}d\lambda_{n}.
\end{split}
\end{equation}

Note that
\begin{equation}
\label{equ:effect20140122.2}
\begin{split}
&(1-e^{-(t_{0}+B_{0})})\int_{D_{v}}\frac{1}{B_{0}}(\mathbb{I}_{\{-t_{0}-B_{0}<t<-t_{0}\}}\circ\psi)|F|^{2}e^{-\psi}d\lambda_{n}
\\&\leq(e^{t_{0}+B_{0}}-1)
\int_{D_{v}}\frac{1}{B_{0}}(\mathbb{I}_{\{-t_{0}-B_{0}<t<-t_{0}\}}\circ\psi)|F|^{2}d\lambda_{n},
\end{split}
\end{equation}
Then we have
\begin{equation}
\label{equ:effect20140122.3}
\begin{split}
&((\int_{D_{v}}|F_{v,t_0}|^{2}d\lambda_{n})^{1/2}-(\int_{D_{v}}|(1-b_{t_0}(\psi))F|^{2}d\lambda_{n})^{1/2})^{2}
\\&\leq(e^{t_{0}+B_{0}}-1)
\int_{D_{v}}\frac{1}{B_{0}}(\mathbb{I}_{\{-t_{0}-B_{0}<t<-t_{0}\}}\circ\psi)|F|^{2}d\lambda_{n},
\end{split}
\end{equation}

Note that
$$e^{\frac{1}{p}t_{0}}\int_{D_{v}}|(1-b_{t_0}(\psi))F|^{2}d\lambda_{n}\leq
e^{\frac{1}{p}t_{0}}\int_{D_{v}}|\mathbb{I}_{\{\frac{1}{p}\psi\leq-\frac{1}{p}t_{0}\}}F|^{2}d\lambda_{n}\leq \int_{D_{v}}|F|^{2}e^{-\frac{1}{p}\psi}d\lambda_{n}$$
As
$$C_{1}\geq\int_{D_{v}}|F|^{2}e^{-\frac{1}{p}\psi}d\lambda_{n}=\int_{D_{v}}|F|^{2}e^{-\varphi}d\lambda_{n}.$$
and $\inf\{\int_{D_{v}}|F_{1}|^{2}d\lambda_{n}|F_{1}\in\mathcal{O}_{D_{v}},(F_{1}-F,z_{0})\in\mathcal{I}(\psi)_{z_{0}}\}\geq C_{2}$,
when
$$e^{\frac{1}{p}t_{0}}\geq \frac{C_{1}}{C_{2}},$$
then we have
\begin{equation}
\label{equ:effect20140122.4}
\begin{split}
&(C_{2}^{1/2}-(C_{1}e^{-\frac{1}{p}t_{0}})^{1/2})^{2}
\\&\leq((\int_{D_{v}}|F_{v,t_0}|^{2}d\lambda_{n})^{1/2}-(\int_{D_{v}}|(1-b_{t_0}(\psi))F|^{2}d\lambda_{n})^{1/2})^{2}.
\end{split}
\end{equation}
It follows that
\begin{equation}
\label{equ:effect20140122.4}
\begin{split}
&(C_{2}^{1/2}-(C_{1}e^{-\frac{1}{p}t_{0}})^{1/2})^{2}
\\&\leq(e^{t_{0}+B_{0}}-1)
\int_{D_{v}}\frac{1}{B_{0}}(\mathbb{I}_{\{-t_{0}-B_{0}<t<-t_{0}\}}\circ\psi)|F|^{2}d\lambda_{n}
\\&\leq e^{t_{0}+B_{0}}\int_{D_{v}}\frac{1}{B_{0}}(\mathbb{I}_{\{-t_{0}-B_{0}<t<-t_{0}\}}\circ\psi)|F|^{2}d\lambda_{n}.
\end{split}
\end{equation}

Replacing $t_{0}$ by $kB_{0}$,
and assuming that $e^{\frac{1}{p}kB_{0}}\geq\frac{C_{1}}{C_{2}}$,
we obtain that
\begin{equation}
\label{equ:effect20140126.1}
\begin{split}
&(C_{2}^{1/2}-(C_{1}e^{-\frac{1}{p}kB_{0}})^{1/2})^{2}
\\&\leq e^{(k+1)B_{0}}\int_{D_{v}}\frac{1}{B_{0}}(\mathbb{I}_{\{-(k+1)B_{0}<t<-kB_{0}\}}\circ\psi)|F|^{2}d\lambda_{n}.
\end{split}
\end{equation}
It follows that
\begin{equation}
\label{equ:effect20140126.3}
\begin{split}
& B_{0}e^{-(k+1)B_{0}}e^{\frac{1}{p}kB_{0}}(C_{2}^{1/2}-(C_{1}e^{-\frac{1}{p}kB_{0}})^{1/2})^{2}
\\&\leq e^{\frac{1}{p}kB_{0}}\int_{D_{v}}(\mathbb{I}_{\{-(k+1)B_{0}<t<-kB_{0}\}}\circ\psi)|F|^{2}d\lambda_{n}.
\end{split}
\end{equation}

Taking $k_{0}$,
such that
$$e^{\frac{1}{p}k_{0}B_{0}}\geq \frac{C_{1}}{C_{2}}\geq e^{\frac{1}{p}(k_{0}-1)B_{0}},$$
and taking sum, we obtain
\begin{equation}
\label{equ:effect20140126.4}
\begin{split}
&\sum_{k=k_{0}}^{+\infty}B_{0}e^{-(k+1)B_{0}}e^{\frac{1}{p}kB_{0}}(C_{2}^{1/2}-(C_{1}e^{-\frac{1}{p}kB_{0}})^{1/2})^{2}
\\&\leq \sum_{k=k_{0}}^{+\infty}e^{\frac{1}{p}kB_{0}}\int_{D_{v}}(\mathbb{I}_{\{-(k+1)B_{0}<t<-kB_{0}\}}\circ\psi)|F|^{2}d\lambda_{n}
\leq\int_{D_{v}}|F|^{2}e^{-\frac{1}{p}\psi}d\lambda_{n}\leq C_{1}.
\end{split}
\end{equation}

Note that
\begin{equation}
\label{equ:effect20140126.5}
\begin{split}
&\sum_{k=k_{0}}^{+\infty}B_{0}e^{-(k+1)B_{0}}e^{\frac{1}{p}kB_{0}}(C_{2}^{1/2}-(C_{1}e^{-\frac{1}{p}kB_{0}})^{1/2})^{2}
\\=&
\sum_{k=k_{0}}^{+\infty}(B_{0}e^{-(k+1)B_{0}}e^{\frac{1}{p}kB_{0}}C_{2}-
2B_{0}e^{-B_{0}}e^{-(\frac{1}{2p}+1)kB_{0}}e^{\frac{1}{p}kB_{0}}C_{2}^{1/2}C_{1}^{1/2}
\\&+B_{0}e^{-B_{0}}e^{-(\frac{1}{p}+1)kB_{0}}e^{\frac{1}{p}kB_{0}}C_{1})
\\=&
\frac{B_{0}}{1-e^{-(1-\frac{1}{p})B_{0}}}e^{-(k_{0}(1-\frac{1}{p})+1)B_{0}}C_{2}
-2\frac{B_{0}}{1-e^{-(1-\frac{1}{2p})B_{0}}}e^{-(k_{0})(1-\frac{1}{2p})B_{0}-B_{0}}C_{2}^{1/2}C_{1}^{1/2}+
\\&+\frac{B_{0}}{1-e^{-B_{0}}}e^{-(k_{0})B_{0}-B_{0}}C_{2}
\\\geq&
\frac{B_{0}}{1-e^{-(1-\frac{1}{p})B_{0}}}(\frac{C_{2}}{C_{1}})^{p(1-\frac{1}{p})}e^{-B_{0}-(1-\frac{1}{p})B_{0}}C_{2}
\\&-2\frac{B_{0}}{1-e^{-(1-\frac{1}{2p})B_{0}}}(\frac{C_{2}}{C_{1}})^{p(1-\frac{1}{2p})}e^{-B_{0}}C_{2}^{1/2}C_{1}^{1/2}
\\&+\frac{B_{0}}{1-e^{-B_{0}}}(\frac{C_{2}}{C_{1}})^{p}e^{-B_{0}-B_{0}}C_{1}
,
\end{split}
\end{equation}

Take limitation
\begin{equation}
\label{equ:effect20140126.6}
\begin{split}
\lim_{B_{0}\to0}&(\frac{B_{0}}{1-e^{-(1-\frac{1}{p})B_{0}}}(\frac{C_{2}}{C_{1}})^{p(1-\frac{1}{p})}e^{-B_{0}-(1-\frac{1}{p})B_{0}}C_{2}
\\&-2\frac{B_{0}}{1-e^{-(1-\frac{1}{2p})B_{0}}}(\frac{C_{2}}{C_{1}})^{p(1-\frac{1}{2p})}e^{-B_{0}}C_{2}^{1/2}C_{1}^{1/2}
\\&+\frac{B_{0}}{1-e^{-B_{0}}}(\frac{C_{2}}{C_{1}})^{p}e^{-B_{0}-B_{0}}C_{1})
\\=&(1-\frac{1}{p})^{-1}(\frac{C_{2}}{C_{1}})^{p(1-\frac{1}{p})}C_{2}-2(1-\frac{1}{2p})^{-1}(\frac{C_{2}}{C_{1}})^{p(1-\frac{1}{2p})}C_{2}^{1/2}C_{1}^{1/2}
+(\frac{C_{2}}{C_{1}})^{p}C_{1}
\\=&(1-\frac{1}{p})^{-1}(\frac{C_{2}}{C_{1}})^{p-1}C_{2}-2(1-\frac{1}{2p})^{-1}(\frac{C_{2}}{C_{1}})^{p-\frac{1}{2}}C_{2}^{1/2}C_{1}^{1/2}
+(\frac{C_{2}}{C_{1}})^{p}C_{1}
\\=&(1-\frac{1}{p})^{-1}(\frac{C_{2}}{C_{1}})^{p}C_{1}-2(1-\frac{1}{2p})^{-1}(\frac{C_{2}}{C_{1}})^{p}C_{1}
+(\frac{C_{2}}{C_{1}})^{p}C_{1}
\\=&((1-\frac{1}{p})^{-1}-2(1-\frac{1}{2p})^{-1}
+1)(\frac{C_{2}}{C_{1}})^{p}C_{1}.
\end{split}
\end{equation}
By inequality \ref{equ:effect20140126.4}, \ref{equ:effect20140126.5} and \ref{equ:effect20140126.6},
it follows that

\begin{equation}
\label{equ:effect20140126.7}
\begin{split}
((1-\frac{1}{p})^{-1}-2(1-\frac{1}{2p})^{-1}
+1)(\frac{C_{2}}{C_{1}})^{p}C_{1}\leq C_{1},
\end{split}
\end{equation}
that is to say
\begin{equation}
\label{equ:effect20140126.8}
\begin{split}
(\frac{1}{(p-1)(2p-1)})^{\frac{1}{p}}\leq \frac{C_{1}}{C_{2}},
\end{split}
\end{equation}
Thus we obtain Proposition \ref{p:effect}.

\subsection{Proof of Theorem \ref{t:effect}}
$\\$

Note that $C_{F,p\varphi}(z_{0})>0$.
It follows that $(F,z_{0})\not\in\mathcal{I}(p\varphi)_{z_{0}}$,
i.e.
$|F|^{2}e^{-p\varphi}$ is not integrable near $z_{0}$.
Note that
\begin{equation}
\begin{split}
C_{F,p\varphi}(z_{0})
&\geq\inf\{\|F_{1}\|^{2}_{0}|(F_{1}-F,z_{0})\in\mathcal{I}_{+}(2c_{z_{0}}^{F}(\varphi)\varphi)_{z_0}{\,}\&{\,}F_{1}\in\mathcal{O}(D)\}>0.
\end{split}
\end{equation}
Using Proposition \ref{p:effect}, we obtain Theorem \ref{t:effect}.

\section{Proof of the  lower semicontinuity property of plurisubharmonic functions with a multiplier}
In this section, we explicitly point out that the proof of Proposition \ref{p:lower}
is implicitly contained in our proof of Proposition 1.6 in \cite{GZopen-b}.

For the sake of completeness, we recall our proof of Proposition 1.6 in \cite{GZopen-b} with the following slightly
modifications:

$1)$ changing $\mu$ into $\mu_{F}$, and $\mu(\Omega)$ into $C_{0}$ respectively;

$2)$ changing $e^{-\phi}$ into $|F|^{2}e^{-\phi}$,
$\\$
where $C_{0}:=\inf_{m}\inf\{\|F_{1}\|^{2}_{0} |(F_{1}-F,z_{0})\in\mathcal{I}(\phi_{m})_{z_0}\&F_{1}\in\mathcal{O}(\Omega)\}>0$,
and $\mu_{F}:=|F|^{2}\mu$ on $\Delta^{n}$ (the choice of $\Omega$ see the following part of the present section).

Let's recall our proof of Proposition 1.6 in \cite{GZopen-b} in details:

We prove Proposition \ref{p:lower} by contradiction.
If $|F|^{2}e^{-\phi}$ is integrable near $o\in\Delta^{n}$,
then there exists a strong pseudoconvex domain $\Omega\subset\subset\Delta^{n}$,
such that $|F|^{2}e^{-\phi}$ is $L^{1}$ integrable on $\Omega$.

Without losing of generality,
we assume that $\Omega=\mathbb{B}(o,r)$,
where $r>0$ small enough.

As $|F|^{2}e^{-\phi}$ is $L^{1}$ integrable on $\Omega$,
then
\begin{equation}
\label{equ:20131121a}\lim_{R\to+\infty}e^{R}\mu_{F}(\{\phi<-R\})=0.
\end{equation}
Therefore there exists $t_{1}>0$,
such that
\begin{equation}
\label{equ:20131121f}
\mu_{F}(\{\phi<-t_{1}+1\})<\frac{1}{6}C_{0}.
\end{equation}

As $\{\phi_{m}\}_{m=1,2,\cdots}$ is convergent to $\phi$,
and $\{F_{m}\}_{m=1,2,\cdots}$ with uniform bound, it follows that there exists $m_{0}>0$,
such that for any $m\geq m_{0}$,
\begin{equation}
\label{equ:20131121g}
\mu_{F_{m}}(\{|\phi_{m}-\phi|\geq 1\})<\frac{1}{12}C_{0}.
\end{equation}

Note that
$$(\{\phi_{m}<-t_{1}\}\setminus\{|\phi_{m}-\phi|\geq 1\})\subset\{\phi<-t_{1}+1\},$$
for any $m\geq m_{0}$.
Therefore
\begin{equation}
\label{equ:20131121b}
\mu_{F_{m}}(\{\phi_{m}<-t_{1}\})\leq\mu_{F_{m}}(\{\phi<-t_{1}+1\})+\mu_{F_{m}}(\{|\phi_{m}-\phi|\geq 1\})<\frac{1}{4}C_{0},
\end{equation}
for any $m\geq m_{0}$.

In Lemma \ref{p:GZ_JM_sharp} with $B_{0}=1$,
then there exists a holomorphic function $F_{v,t_{0}}$ on $\Omega$,
satisfying:
\begin{equation}
\label{equ:20131121c}
(F_{m,v,t_{0}}-F_{m},o)\in\mathcal{I}(\phi_{m})_{o}
\end{equation}
and
\begin{equation}
\label{equ:20131026b}
\begin{split}
&\int_{\Omega}|F_{m,v,t_0}-(1-b_{t_0}(\phi_{m}))F|^{2}d\lambda_{n}
\\\leq&\int_{\Omega}(\mathbb{I}_{\{-t_{0}-1< t<-t_{0}\}}\circ\phi_{m})|F_{m}|^{2}e^{-\phi_{m}}d\lambda_{n}.
\end{split}
\end{equation}

It follows from \ref{equ:20131121c}
that
\begin{equation}
\label{equ:20131027a}
\begin{split}
\int_{\Omega}|F_{m,v,t_0}|^{2}d\lambda_{n}\geq C_{0}.
\end{split}
\end{equation}

It follows from inequality \ref{equ:20131121b}, equality \ref{equ:20131121c},
and $b_{t_{0}}(t)|_{\{t\geq -t_{0}\}}=1$,
that
\begin{equation}
\label{equ:20131121d}
\begin{split}
\int_{\Omega}|(1-b_{t_0}(\phi_{m}))F_{m}|^{2}d\lambda_{n}
&=\int_{\{\phi<-t_{0}\}\cap\Omega}|(1-b_{t_0}(\phi_{m}))F_{m}|^{2}d\lambda_{n}
\\&\leq\int_{\{\phi<-t_{0}\}\cap\Omega}|F_{m}|^{2}d\lambda_{n}
\\&\leq\mu_{F_{m}}(\{\phi<-t_{0}\})<\frac{1}{4}C_{0}
\end{split}
\end{equation}

It follows from
inequalities \ref{equ:20131121d} and \ref{equ:20131027a}
that
\begin{equation}
\label{equ:20131026a}
\begin{split}
&(\int_{\Omega}|F_{m,v,t_0}-(1-b_{t_0}(\phi_{m}))F_{m}|^{2}d\lambda_{n})^{1/2}
\\&\geq(\int_{\Omega}|F_{m,v,t_0}|^{2}d\lambda_{n})^{1/2}-
(\int_{\Omega}|(1-b_{t_0}(\phi_{m}))F_{m}|^{2}d\lambda_{n})^{1/2}
\\&\geq C_{0}^{1/2}-2^{-1}C_{0}^{1/2}=2^{-1}C_{0}^{1/2},
\end{split}
\end{equation}
for any $t_{0}>t_{1}$ and any $m\geq m_{0}$.

It follows from inequalities \ref{equ:20131026b} and \ref{equ:20131026a}
that
$$\int_{\Omega}(\mathbb{I}_{\{-t_{0}-1< t<-t_{0}\}}\circ\phi_{m})|F_{m}|^{2}e^{-\phi_{m}}d\lambda_{n}\geq 2^{-2}C_{0},$$
for any $t_{0}>t_{1}$ and any $m\geq m_{0}$.

Note that
$$\mu_{F_{m}}({\{-t_{0}-1<\phi_{m}<-t_{0}\}})e^{t_{0}+1}\geq \int_{\Omega}(\mathbb{I}_{\{-t_{0}-1< t<-t_{0}\}}\circ\phi_{m})|F_{m}|^{2}e^{-\phi_{m}}d\lambda_{n},$$
for any $t_{0}>t_{1}$ and any $m\geq m_{0}$.
Therefore
$$\mu_{F_{m}}({\{-t_{0}-1<\phi_{m}<-t_{0}\}})\geq e^{-t_{0}-1}2^{-2}C_{0},$$
for any $t_{0}>t_{1}$ and any $m\geq m_{0}$.

As $\{\phi_{m}\}_{m=1,2,\cdots}$ is convergent to $\phi$ in Lebesgue measure,
and  $\{F_{m}\}_{m=1,2,\cdots}$ is uniformly bounded,
then
there exists large enough positive integer $m_{1}\geq m_{0}$,
such that
$$\mu_{F_{m_{1}}}(\{|\phi_{m_1}-\phi|\geq1\})<\frac{1}{2}e^{-t_{0}-1}2^{-2}C_{0},$$
for any $t_{0}>t_{1}$.

Note that
$$\{\phi<-t_{0}+1\}\supset(\{-t_{0}-1<\phi_{m_{1}}<-t_{0}\}\setminus\{|\phi_{m_1}-\phi|\geq1\}).$$
Then we have
\begin{equation}
\begin{split}
\mu_{F_{m_{1}}}({\{\phi<-t_{0}+1\}})&\geq\mu_{F_{m_{1}}}(\{-t_{0}-1<\phi_{m_{1}}<-t_{0}\}\setminus\{|\phi_{m_1}-\phi|\geq1\})
\\&\geq \mu_{F_{m_{1}}}(\{-t_{0}-1<\phi_{m_{1}}<-t_{0}\})-\mu_{F_{m_{1}}}(\{|\phi_{m_1}-\phi|\geq1\})
\\&\geq\frac{1}{2}e^{-t_{0}-1}2^{-2}C_{0},
\end{split}
\end{equation}
for any $t_{0}>t_{1}$,
i.e.
\begin{equation}
\label{equ:20131121e}
e^{t_{0}-1}\mu_{F_{m_{1}}}({\{\phi<-t_{0}+1\}})\geq\frac{1}{2}e^{-2}2^{-2}C_{0},
\end{equation}
for any $t_{0}>t_{1}$.

When $m_{1}$ goes to infinity,
as $\{F_{m}\}_{m=1,2,\cdots}$ is convergent to $F$ in Lebesgue measure with uniform bound,
by the dominated convergence theorem,
it follows that
$$\lim_{m_{1}\to+\infty}e^{t_{0}-1}\mu_{F_{m_{1}}}({\{\phi<-t_{0}+1\}})=e^{t_{0}-1}\mu_{F}({\{\phi<-t_{0}+1\}}).$$

Using inequality \ref{equ:20131121e}, we obtain
$$e^{t_{0}-1}\mu_{F}({\{\phi<-t_{0}+1\}})\geq\frac{1}{2}e^{-2}2^{-2}C_{0},$$
which contradicts to equality \ref{equ:20131121a}.

Proposition \ref{p:lower} has thus been proved.

\begin{Remark}
Using the same method as in the above proof with more subtle bounds
in inequalities \ref{equ:20131121f} and \ref{equ:20131121g}, one can
obtain inequality \ref{equ:20131121e} with a lower bound
$e^{-2}C_{0}$.
\end{Remark}

\section{Proofs of the effectiveness of Conjecture D-K and Conjecture J-M}

In this section, we explicitly point out that
Proposition \ref{p:effect_conj.D-K} and Proposition \ref{p:effect_conj.J-M}
are already implicitly contained in \cite{GZopen-b}.

\subsection{Proof of Proposition \ref{p:effect_conj.D-K}}
$\\$

In this subsection, we explicitly point out that the proof of Proposition \ref{p:effect_conj.D-K}
is implicitly contained in our proof of Theorem \ref{t:GZ_JM} in \cite{GZopen-b}.

For the sake of completeness, we recall some steps of our proof of Theorem \ref{t:GZ_JM} in \cite{GZopen-b} with slightly
modification: changing $o$ into $z_{0}$.

By Proposition \ref{p:GZ_JM_sharp},
it follows that there exists $F_{v,t_{0}}$,
which is a holomorphic function on $D_{v}$
satisfying:

\begin{equation}
\label{equ:13.10.11a}
\begin{split}
&\int_{ D_v}|F_{v,t_0}-(1-b_{t_0}(\varphi))F|^{2}d\lambda_{n}
\\\leq&\int_{D_v}\mathbb{I}_{\{-t_{0}-B_{0}<\varphi<-t_{0}\}}|F|^{2}e^{-\varphi}d\lambda_{n}.
\end{split}
\end{equation}
and
$$(F_{v,t_{0}}-F,z_{0})\in\mathcal{I}(\varphi)_{z_{0}}.$$

Note that
\begin{equation}
\label{equ:20131120b}
\begin{split}
(\int_{ D_v}|F_{v,t_0}|^{2}d\lambda_{n})^{\frac{1}{2}}&\leq
(\int_{ D_v}|F_{v,t_0}-(1-b_{t_0}(\varphi))F|^{2}d\lambda_{n})^{\frac{1}{2}}
\\&+(\int_{ D_v}|(1-b_{t_0}(\varphi))F|^{2}d\lambda_{n})^{\frac{1}{2}},
\end{split}
\end{equation}
and
\begin{equation}
\label{equ:20140227a}
\begin{split}
\lim_{t_{0}\to0}\int_{ D_v}|(1-b_{t_0}(\varphi))F|^{2}d\lambda_{n}=0,
\end{split}
\end{equation}
then it follows that
\begin{equation}
\label{equ:20140227b}
\begin{split}
&\liminf_{t_{0}\to+\infty}\int_{D_v}\mathbb{I}_{\{-t_{0}-B_{0}<\varphi<-t_{0}\}}|F|^{2}e^{t_{0}+B_{0}}d\lambda_{n}
\\&\geq\liminf_{t_{0}\to+\infty}\int_{D_v}\mathbb{I}_{\{-t_{0}-B_{0}<\varphi<-t_{0}\}}|F|^{2}e^{-\varphi}d\lambda_{n}
\\&\geq\inf\{\parallel F_{1}\parallel_{0}^{2}|(F_{1}-F,z_{0})\in\mathcal{I}(\varphi)_{z_{0}}\,\&\,F_{1}\in\mathcal{O}_{D_{v}}\}
\end{split}
\end{equation}

Proposition \ref{p:effect_conj.D-K} has thus been proved.

\subsection{A proposition used in the proof of Conjecture J-M}
$\\$

Let
$$\varphi:=(1+\delta)\max\{\psi,\log|F|^{2}\},$$
and
$$\Psi:=\min\{\psi-\log|F|^{2},0\},$$
where $\delta$ is a positive integer.
Then $\Psi+\varphi$ and $(1+\delta)\Psi+\varphi$ are both plurisubharmonic functions on $\Delta^{n}$.

In \cite{GZopen-b}, we proved Theorem \ref{t:GZ_JM201312} by the following proposition:

\begin{Proposition} \label{p:GZ_JM201312}\cite{GZopen-b}
Let $D_{v}$ be a strongly pseudoconvex domain relatively compact in $\Delta^{n}$ containing $o$.

Then there exists  a holomorphic function $F_{v,t_{0}}$ on $D_{v}$,
satisfying:
\begin{equation}
\label{equ:20131210a}
(F_{v,t_{0}}-F^{1+\delta},o)\in\mathcal{I}(\varphi+\Psi)_{o}
\end{equation}
and
\begin{equation}
\label{equ:3.4}
\begin{split}
&\int_{ D_v}|F_{v,t_0}-(1-b_{t_0}(\Psi))F^{1+\delta}|^{2}e^{-\varphi}d\lambda_{n}
\\\leq&(1+\frac{1}{\delta})\int_{D_v}(\frac{1}{B_{0}}\mathbb{I}_{\{-t_{0}-B_{0}< t<-t_{0}\}}\circ\Psi)|F^{1+\delta}|^{2}e^{-\varphi}e^{-\Psi}d\lambda_{n},
\end{split}
\end{equation}
where
$$b_{t_{0}}:=\int_{-\infty}^{t}\frac{1}{B_{0}}\mathbb{I}_{(-t_{0}-B_{0},-t_{0})}ds,$$
and $t_{0}\geq0.$
\end{Proposition}

Replacing the strong pseudoconvexity of $D_{v}$ by pseudoconvexity, and $o$ by $z_{0}\in D_{v}$ for any $v$,
Lemma \ref{p:GZ_JM_sharp} also holds.

\subsection{Proof of Proposition \ref{p:effect_conj.J-M}}
$\\$

In this subsection, we explicitly point out that the proof of Proposition \ref{p:effect_conj.D-K}
is implicitly contained in our proof of Theorem \ref{t:GZ_JM} in \cite{GZopen-b}.

For the sake of completeness, we recall our proof of Theorem \ref{t:GZ_JM201312} in \cite{GZopen-b} with slightly
modification: changing $o$ into $z_{0}$.

By Proposition \ref{p:GZ_JM201312} with $\Psi:=\min\{\psi-\log|F|^{2},0\}$,
it follows that there exists $F_{v,t_{0}}$,
which is a holomorphic function on $D_{v}$
satisfying:

\begin{equation}
\label{equ:13.10.11a}
\begin{split}
&\int_{ D_v}|F_{v,t_0}-(1-b_{t_0}(\Psi))F^{1+\delta}|^{2}e^{-\varphi}d\lambda_{n}
\\\leq&(1+\frac{1}{\delta})\int_{D_v}\frac{1}{B_{0}}\mathbb{I}_{\{-t_{0}-B_{0}<\Psi<-t_{0}\}}|F^{1+\delta}|^{2}e^{-\varphi-\Psi}d\lambda_{n},
\end{split}
\end{equation}
and
$$(F_{v,t_{0}}-F^{1+\delta},z_{0})\in\mathcal{I}(\varphi+\Psi)_{z_{0}}.$$

Note that
\begin{equation}
\label{equ:20131120b}
\begin{split}
(\int_{ D_v}|F_{v,t_0}|^{2}e^{-\varphi}d\lambda_{n})^{\frac{1}{2}}&\leq
(\int_{ D_v}|F_{v,t_0}-(1-b_{t_0}(\Psi))F^{1+\delta}|^{2}e^{-\varphi}d\lambda_{n})^{\frac{1}{2}}\\&+(\int_{ D_v}|(1-b_{t_0}(\Psi))F^{1+\delta}|^{2}e^{-\varphi}d\lambda_{n})^{\frac{1}{2}},
\end{split}
\end{equation}
and
\begin{equation}
\label{equ:20140227c}
\begin{split}
\lim_{t_{0}\to0}\int_{ D_v}|(1-b_{t_0}(\Psi))F^{1+\delta}|^{2}e^{-\varphi}d\lambda_{n}=0,
\end{split}
\end{equation}
then it follows that
\begin{equation}
\label{equ:20140227d}
\begin{split}
&(1+\frac{1}{\delta})e^{B_{0}+1}\frac{1}{B_{0}}\liminf_{t_{0}\to\infty}\int_{D_v}\mathbb{I}_{\{-t_{0}-B_{0}<\Psi<-t_{0}\}}e^{t_{0}}d\lambda_{n}
\\&\geq\liminf_{t_{0}\to\infty}(1+\frac{1}{\delta})
\int_{D_v}\frac{1}{B_{0}}\mathbb{I}_{\{-t_{0}-B_{0}<\Psi<-t_{0}\}}|F^{1+\delta}|^{2}e^{-\varphi-\Psi}d\lambda_{n}
\\&\geq\frac{1}{e^{\sup_{D}\varphi}}\inf\{\parallel F_{1}\parallel_{0}^{2}|(F_{1}-F^{1+\delta},z_{0})\in\mathcal{I}(\varphi+\Psi)_{z_{0}}\,\&\,F_{1}\in\mathcal{O}_{D_{v}}\}.
\end{split}
\end{equation}

Note that $\varphi+\Psi=\psi+\delta\max\{\psi,\log|F|^{2}\}$.
Proposition \ref{p:effect_conj.J-M} has thus been proved.

\subsection{Proof of Remark \ref{r:JM-1}}
$\\$

When $D=\Delta$, $F=1$, $\psi=\log|z|^{2}$,
then $I$ is trivial, $\log|I|=0$ and $c_{o}^{I}(\log|z|^{2})=c_{o}(\log|z|^{2})=\frac{1}{2}$.

As $f_{i}=1$, then it follows that $\max\{2c_{o}^{I}(\psi)\psi,2\log|f_{i}|\}=0$ for any $i$,
and $\max\{2c_{o}^{I}(\psi)\psi,2\log|I|\}=0.$

It is clear that
$$K_{2c_{o}^{I}(\psi)\psi+\delta\max\{2c_{o}^{I}(\psi)\psi,\log|f_{i}|^{2}\},F^{1+\delta}}(o)=K_{\psi,1}(o)=K(o)=\frac{1}{\pi},$$
and
$$\liminf_{r\to0}\frac{1}{r^2}\mu(\{c^{I}_{o}(\psi)\psi-\log|I|<\log r\})=\lim_{r\to0}\frac{1}{r^2}\mu(\{\log|z|<\log r\})=\frac{1}{\pi}.$$

When $\delta$ goes to $\infty$, the equality in inequality \ref{equ:20140304} holds.

\subsection{Proof of Remark \ref{r:JM-DK}}
$\\$

When $F=1$, $\psi<0$, then $f_{i}=1$ and $c_{o}^{I}(\psi)=c_{o}(\psi)$.

It is clear that
$$K^{-1}_{2c_{o}^{I}(\psi)\psi+\delta\max\{2c_{o}^{I}(\psi)\psi,\log|f_{i}|^{2}\},1}(o)
=K^{-1}_{2c_{o}(\psi)\psi,1}(o)\geq K^{-1}(o),$$
and
$$\sup_{D}e^{(1+\delta)\max\{2c_{o}^{I}(\psi)\psi,2\log|I|\}}=\sup_{D}e^{(1+\delta)\max\{2c_{o}^{I}(\psi)\psi,0\}}=\sup_{D}e^{0}=1.$$

Replacing $r$ in inequality \ref{equ:20140304} by $r^{2c_{o}(\psi)}$,
the optimal effectiveness of conjecture J-M degenerates to the optimal effectiveness of conjecture D-K:

\begin{equation}
\begin{split}
\liminf_{r\to0}\frac{1}{r^{2c_{o}(\psi)}}\mu(\{\psi<\log r\})=&\liminf_{r\to0}\frac{1}{r^2}\mu(\{c^{I}_{o}(\psi)\psi-\log|I|<\log r\})
\\&\geq \sup_{\delta\in\{1,2,\cdots\}}\frac{K^{-1}(o)}{1+\frac{1}{\delta}}=K^{-1}(o),
\end{split}
\end{equation}

\section{Proofs of preparatory results}

In this section, we recall some main steps in our proof in
\cite{guan-zhou13p} (see also \cite{guan-zhou13ap,GZopen-b}) with some slight
modifications in order to prove Lemma \ref{p:GZ_JM_sharp}.

\subsection{Proof of Lemma \ref{p:GZ_JM_sharp}}
$\\$

For the sake of completeness, let's recall some steps in our proof in
\cite{guan-zhou13p} (see also \cite{guan-zhou13ap,GZopen-b}) with some slight
modifications in order to prove Lemma \ref{p:GZ_JM_sharp}.

Let $\{v_{t_0,\varepsilon}\}_{t_{0}\in\mathbb{R},\varepsilon\in(0,\frac{1}{8}B_{0})}$ be a family of smooth increasing convex functions on $\mathbb{R}$,
which are continuous functions on $\mathbb{R}\cup\{-\infty\}$, such that:

 $1).$ $v_{t_{0},\varepsilon}(t)=t$ for $t\geq-t_{0}-\varepsilon$, $v_{t_{0},\varepsilon}(t)=constant$ for $t<-t_{0}-B_{0}+\varepsilon$;

 $2).$ $v''_{t_0,\varepsilon}(t)$ are pointwise convergent to $\frac{1}{B_{0}}\mathbb{I}_{(-t_{0}-B_{0},-t_{0})}$, when $\varepsilon\to 0$, and $0\leq v''_{t_0,\varepsilon}(t)\leq 2$ for any $t\in \mathbb{R}$;

 $3).$ $v'_{t_0,\varepsilon}(t)$ are pointwise convergent to $b_{t_{0}}(t)=\int_{-\infty}^{t}\frac{1}{B_{0}}\mathbb{I}_{(-t_{0}-B_{0},-t_{0})}ds$ ($b_{t_{0}}$ is also a continuous function on $\mathbb{R}\cup\{-\infty\}$), when $\varepsilon\to 0$, and $0\leq v'_{t_0,\varepsilon}(t)\leq1$ for any $t\in \mathbb{R}$.

One can construct the family $\{v_{t_0,\varepsilon}\}_{t_{0}\in\mathbb{R},\varepsilon\in(0,\frac{1}{8})}$ by the setting
\begin{equation}
\label{equ:20140101}
\begin{split}
v_{t_0,\varepsilon}(t):=&\int_{-\infty}^{t}(\int_{-\infty}^{t_{1}}(\frac{1}{1-4\varepsilon}
\frac{1}{B_{0}}\mathbb{I}_{(-t_{0}-B_{0}+2\varepsilon,-t_{0}-2\varepsilon)}*\rho_{\frac{1}{4}\varepsilon})(s)ds)dt_{1}
\\&-\int_{-\infty}^{0}(\int_{-\infty}^{t_{1}}(\frac{1}{1-4\varepsilon}\frac{1}{B_{0}}\mathbb{I}_{(-t_{0}-B_{0}+2\varepsilon,
-t_{0}-2\varepsilon)}*\rho_{\frac{1}{4}\varepsilon})(s)ds)dt_{1},
\end{split}
\end{equation}
where $\rho_{\frac{1}{4}\varepsilon}$ is the kernel of convolution satisfying $supp(\rho_{\frac{1}{4}\varepsilon})\subset (-\frac{1}{4}\varepsilon,\frac{1}{4}\varepsilon)$.
Then it follows that
$$v''_{t_0,\varepsilon}(t)=\frac{1}{1-4\varepsilon}\frac{1}{B_{0}}\mathbb{I}_{(-t_{0}-B_{0}+2\varepsilon,-t_{0}-2\varepsilon)}*\rho_{\frac{1}{4}\varepsilon}(t),$$
and
$$v'_{t_0,\varepsilon}(t)=\int_{-\infty}^{t}(\frac{1}{1-4\varepsilon}\frac{1}{B_{0}}\mathbb{I}_{(-t_{0}-B_{0}+2\varepsilon,-t_{0}-2\varepsilon)}
*\rho_{\frac{1}{4}\varepsilon})(s)ds.$$

As $D_{v}\subset\subset\Delta^{n}\subset \mathbb{C}^{n}$,
then there exist negative smooth plurisubharmonic functions $\{\psi_{m}\}_{m=1,2,\cdots}$ on a neighborhood of $\overline{D}_{v}$,
such that the sequence $\{\psi_{m}\}_{m=1,2,\cdots}$ is decreasingly convergent to $\psi$
on a smaller neighborhood of $\overline{D}_{v}$,
when $m\to+\infty$.

Take undetermined functions $s$ and $u$ which will be naturally led to an ODE system after calculations based on the twisted Bochner-Kodaira identity and a lemma of Berndtsson's, and will be determined by solving the ODE system (The idea goes back to \cite{guan-zhou-zhu10}, \cite{guan-zhou12}).

Let $\eta=s(-v_{t_{0},\varepsilon}\circ\psi_{m})$ and $\phi=u(-v_{t_{0},\varepsilon}\circ\psi_{m})$,
where $s\in C^{\infty}((0,+\infty))$ satisfies $s\geq0$, and
$u\in C^{\infty}((0,+\infty))$ satisfies $\lim_{t\to+\infty}u(t)=0$, such that $u''s-s''>0$, and $s'-u's=1$.

Let $\Phi=\psi_{m}+\phi$.

Now let $\alpha=\sum^{n}_{j=1}\alpha_{j}d\bar z^{j}\in Dom_{D_v}
(\bar\partial^*)\cap Ker(\bar\partial)\cap C^\infty_{(0,1)}(\overline {D_v})$.
By Cauchy-Schwarz inequality, it follows that
\begin{equation}
\label{equ:20131130a}
\begin{split}
2\mathrm{Re}(\bar\partial^*_\Phi\alpha,\alpha\llcorner(\bar\partial\eta)^\sharp )_{\Omega,\Phi}
\geq
&-\int_{D_v}g^{-1}|\bar{\partial}^{*}_{\Phi}\alpha|^{2}e^{-\Phi}d\lambda_{n}
\\&+
\sum_{j,k=1}^{n}\int_{D_v}
(-g(\partial_{j} \eta)\bar\partial_{k} \eta )\alpha_{\bar{j} }\overline{{\alpha}_{\bar {k}}}e^{-\Phi}d\lambda_{n}.
\end{split}
\end{equation}

Using the twisted Bochner-Kodaira identity (see Lemma 3.2. in \cite{guan-zhou-zhu10}) and inequality \ref{equ:20131130a},
since $s\geq0$ and $\psi_{m}$ is a plurisubharmonic function on $\overline{D}_{v}$,
we get

\begin{equation}
\label{equ:4.1}
\begin{split}
&\int_{D_v}(\eta+g^{-1})|\bar{\partial}^{*}_{\Phi}\alpha|^{2}e^{-\Phi}d\lambda_{n}
\\\geq&\sum_{j,k=1}^{n}\int_{D_v}
(-\partial_{j}\bar{\partial}_{k}\eta+\eta\partial_{j}\bar{\partial}_{k}\phi-g(\partial_{j} \eta)\bar\partial_{k} \eta )\alpha_{\bar{j} }\overline{{\alpha}_{\bar {k}}}e^{-\Phi}d\lambda_{n},
\end{split}
\end{equation}
where $g$ is a positive continuous function on $D_{v}$.
We need some calculations to determine $g$.

We have

\begin{equation}
\label{}
\begin{split}
&\sum_{1\leq j,k\leq n}(-\partial_{j}\bar{\partial}_{k}\eta+\eta\partial_{j}\bar{\partial}_{k}\phi-g(\partial_{j} \eta)
\bar\partial_{k} \eta )\alpha_{\bar{j} }\overline{{\alpha}_{\bar{k}}}
\\=&(s'-su')\sum_{1\leq j,k\leq n}((v'_{t_0,\varepsilon}\circ\psi_{m})\partial_{j}\bar{\partial}_{k}\psi_{m}+(v''_{t_0,\varepsilon}\circ \psi_{m})\partial_{j}(\psi_{m})\bar{\partial}_{k}(\psi_{m}))\alpha_{\bar{j} }\overline{{\alpha}_{\bar{k}}}
\\+&((u''s-s'')-gs'^{2})\sum_{1\leq j,k\leq n}\partial_{j}
(-v_{t_0,\varepsilon}\circ \psi_{m})\bar{\partial}_{k}(-v_{t_0,\varepsilon}\circ \psi_{m})\alpha_{\bar{j} }\overline{{\alpha}_{\bar{k}}}.
\end{split}
\end{equation}
We omit composite item $(-v_{t_0,\varepsilon}\circ \psi_{m})$ after $s'-su'$ and $(u''s-s'')-gs'^{2}$ in the above equalities.

Set $$g=\frac{u''s-s''}{s'^{2}}\circ(-v_{t_0,\varepsilon}\circ \psi_{m}).$$
It follows that $$\eta+g^{-1}=(s+\frac{s'^{2}}{u''s-s''})\circ(-v_{t_0,\varepsilon}\circ \psi_{m}).$$

Because of $v'_{t_0,\varepsilon}\geq 0$  and $s'-su'=1$, using inequalities \ref{equ:4.1}, we have
\begin{equation}
\label{equ:3.1}
\begin{split}
\int_{D_v}(\eta+g^{-1})|\bar{\partial}^{*}_{\Phi}\alpha|^{2}e^{-\Phi}d\lambda_{n}
\geq\int_{D_v}(v''_{t_0,\varepsilon}\circ{\psi_{m}})
\big|\alpha\llcorner(\bar \partial \psi_{m})^\sharp\big|^2e^{-\Phi}d\lambda_{n}.
\end{split}
\end{equation}

Let $\lambda=\bar{\partial}[(1-v'_{t_0,\varepsilon}(\psi_{m})){F}]$.
By the definition of contraction, Cauchy-Schwarz inequality and inequality \ref{equ:3.1},
it follows that
\begin{equation}
\begin{split}
&|(\lambda,\alpha)_{D_v,\Phi}|^{2}=|((v''_{t_0,\varepsilon}\circ{\psi_{m}})\bar\partial\psi_{m} F,\alpha)_{D_v,\Phi}|^{2}
\\=&|((v''_{t_0,\varepsilon}\circ{\psi_{m}})F,\alpha\llcorner(\bar\partial\psi_{m})^\sharp\big)_{D_v,\Phi}|^{2}
\\\leq&\int_{D_v}
(v''_{t_0,\varepsilon}\circ{\psi_{m}})| F|^2e^{-\Phi}d\lambda_{n}\int_{D_v}(v''_{t_0,\varepsilon}\circ{\psi_{m}})
\big|\alpha\llcorner(\bar\partial\psi_{m})^\sharp\big|^2e^{-\Phi}d\lambda_{n}
\\\leq&
(\int_{D_v}
(v''_{t_0,\varepsilon}\circ{\psi_{m}})| F|^2e^{-\Phi}d\lambda_{n})
(\int_{ D_v}(\eta+g^{-1})|\bar{\partial}^{*}_{\Phi}\alpha|^{2}e^{-\Phi}d\lambda_{n}).
\end{split}
\end{equation}

Let $\mu:=(\eta+g^{-1})^{-1}$. Using a lemma of Berndtsson's (see \cite{berndtsson} or Lemma 3.7. in \cite{guan-zhou-zhu10}),
we have locally $L^{1}$ function $u_{v,t_0,m,\varepsilon}$ on $D_{v}$ such that $\bar{\partial}u_{v,t_0,m,\varepsilon}=\lambda$,
and
\begin{equation}
 \label{equ:3.2}
 \begin{split}
 &\int_{ D_v}|u_{v,t_0,m,\varepsilon}|^{2}(\eta+g^{-1})^{-1} e^{-\Phi}d\lambda_{n}
  \leq\int_{D_v}(v''_{t_0,\varepsilon}\circ{\psi_{m}})| F|^2e^{-\Phi}d\lambda_{n}.
  \end{split}
\end{equation}

Let $\mu_{1}=e^{v_{t_0,\varepsilon}\circ\psi_{m}}$, $\tilde{\mu}=\mu_{1}e^{\phi}$.
Assume that we can choose $\eta$ and $\phi$ such that $\tilde{\mu}\leq \mathbf{C}(\eta+g^{-1})^{-1}=\mu$, where $\mathbf{C}=1$.

Note that $v_{t_0,\varepsilon}(\psi_{m})\geq\psi_{m}$.
Then it follows that
\begin{equation}
\label{equ:3.8}
\begin{split}
\int_{ D_v}|u_{v,t_0,m,\varepsilon}|^{2} d\lambda_{n}\leq\int_{ D_v}|u_{v,t_0,m,\varepsilon}|^{2}\mu_{1}e^{\phi} e^{-\psi_{m}-\phi}d\lambda_{n}
=\int_{ D_v}|u_{v,t_0,m,\varepsilon}|^{2}\tilde{\mu}e^{-\Phi}d\lambda_{n}.
\end{split}
\end{equation}

Using inequalities \ref{equ:3.2} and \ref{equ:3.8},
we obtain that
$$\int_{ D_v}|u_{v,t_0,m,\varepsilon}|^{2} d\lambda_{n}\leq\mathbf{C}\int_{D_v}
(v''_{t_0,\varepsilon}\circ\psi_{m})| F|^2e^{-\Phi}d\lambda_{n},$$
under the assumption $\tilde{\mu}\leq\mathbf{C} (\eta+g^{-1})^{-1}$.

As $-v_{t_0,\varepsilon}\circ{\psi_{m}}(\overline{D_{v}})\subset\subset(0,t_{0}+1)$
and $\{\psi_{m}\}_{m=1,2,\cdots}$ is decreasing,
then it is clear that
\begin{equation}
\label{equ:20131130b}
-v_{t_0,\varepsilon}\circ{\psi_{m}}(\overline{D_{v}})\subset\subset K_{t_{0}}\subset\subset(0,t_{0}+1)
\end{equation}
where $K_{t_{0}}$ is independent of $m$ and $\varepsilon\in(0,\frac{1}{8}B_{0})$.

As $u$ is positive and smooth on $(0,+\infty)$,
it follows that $\phi$ is uniformly bounded on $\overline{D}_{v}$ independent of $m$.

As $Supp(v''_{t_{0},\varepsilon})\subset\subset(-t_{0}-B_{0},-t_{0})$,
then it is clear that $(v''_{t_0,\varepsilon}\circ{\psi_{m}})| F|^2e^{-\psi_{m}}$ are uniformly bounded on $\overline{D}_{v}$ independent of $m$.

Therefore
$$\int_{D_v}
(v''_{t_0,\varepsilon}\circ{\psi_{m}})| F|^2e^{-\Phi}d\lambda_{n}$$ are uniformly bounded independent of $m$,
for any given $v$, $t_0$, $\varepsilon$.

By weakly compactness of the unit ball of $L^{2}(D_{v})$ and dominated convergence theorem, when $m\to+\infty$,
it follows that the weak limit of some weakly convergent subsequence of
$\{u_{v,t_{0},m,\varepsilon}\}_{m}$ gives $u_{v,t_0,\varepsilon}$ on $D_v$ satisfying
\begin{equation}
\label{equ:3.3}\int_{ D_v}|u_{v,t_0,\varepsilon}|^{2}d\lambda_{n}\leq\frac{\mathbf{C}}{e^{A_{t_0}}}\int_{D_v}
(v''_{t_0,\varepsilon}\circ{\psi})| F|^2e^{-\psi}d\lambda_{n},
\end{equation}
where $A_{t_0}:=\inf_{t\in(t_{0},t_{0}+B_{0})}\{u(t)\}$.

As $\psi_{m}$ is decreasingly convergent to $\psi$ on $\Delta^{n}$,
and $\psi(o)=-\infty$,
then for any given $t_{0}$ there exists $m_{0}$ and a neighbourhood $U_{0}$ of $o\in D_{v}$ on $\Delta^{n}$,
such that for any $m\geq m_{0}$ and $\varepsilon\in(0,\frac{1}{8}B_{0})$,
$v''_{t_0,\varepsilon}\circ\psi_{m}|_{U_{0}}=0$.

It follows that
$$\bar\partial u_{v,t_0,m,\varepsilon}|_{U_0}=\lambda|_{U_0}=\bar{\partial}[(1-v'_{t_0,\varepsilon}(\psi_{m})){F}]|_{U_0}=-(v''_{t_0,\varepsilon}\circ\psi_{m})F\bar{\partial}\psi_{m}|_{U_{0}}=0.$$
That is to say $u_{v,t_0,m,\varepsilon}|_{U_0}$ are all holomorphic.
Therefore $u_{v,t_0,\varepsilon}|_{U_0}$ is holomorphic.

Recall that the integrals $\int_{ D_v}|u_{v,t_{0},m,\varepsilon}|^{2} d\lambda_{n}$ have a uniform bound independent of $m$,
then we can choose a subsequence with respect to $m$ from the chosen weakly convergent subsequence of $u_{v,t_{0},m,\varepsilon}$,
such that the subsequence is uniformly convergent on any compact subset of $U_0$,
and we still denote the subsequence by $u_{v,t_{0},m,\varepsilon}$ without ambiguity.

By the above arguments,
it follows that the right hand side of inequality \ref{equ:3.2} are uniformly bounded independent of $m$ and $\varepsilon\in(0,\frac{1}{8}B_{0})$.

By inequality \ref{equ:3.2},
it follows that
$$\int_{ D_v}|u_{v,t_0,m,\varepsilon}|^{2}(\eta+g^{-1})^{-1} e^{-\phi-\psi_{m}}d\lambda_{n}$$
are uniformly bounded independent of $m$ and $\varepsilon\in(0,\frac{1}{8}B_{0})$.

Using inequality \ref{equ:20131130b},
we obtain that
$$(\eta+g^{-1})^{-1}=(s(-v_{t_{0},\varepsilon}\circ\psi_{m})+\frac{s'^{2}}{u''s-s''}\circ(-v_{t_0,\varepsilon}\circ\psi_{m}))^{-1}$$
and
$e^{-\phi}=e^{-u(-v_{t_{0},\varepsilon}\circ\psi_{m})}$ have positive uniform lower bounds independent of $m$ and $\varepsilon\in(0,\frac{1}{8}B_{0})$.

Then the integrals
$$\int_{ K_0}|u_{v,t_0,m,\varepsilon}|^{2}e^{-\psi_{m}}d\lambda_{n}$$
have a uniform upper bound independent of $m$ and $\varepsilon\in(0,\frac{1}{8}B_{0})$
for any given compact set $K_{0}\subset\subset U_0\cap D_v$,
and
$$\bar\partial u_{v,t_0,m,\varepsilon}|_{U_0}=0.$$

As
$\psi_{m'}\leq\psi_{m}$ where $m'\geq m$,
it follows that
$$|u_{v,t_0,m',\varepsilon}|^{2}e^{-\psi_{m}}\leq|u_{v,t_0,m',\varepsilon}|^{2}e^{-\psi_{m'}}.$$
Then for any given compact set $K_{0}\subset\subset U_0\cap D_v$,
$\int_{ K_0}|u_{v,t_0,m',\varepsilon}|^{2}e^{-\psi_{m}}d\lambda_{n}$
have a uniform bound independent of $m$ and $m'$.

It is clear that
for any given compact set $K_{0}\subset\subset U_0\cap D_v$,
the integrals
$$\int_{ K_0}|u_{v,t_0,\varepsilon}|^{2}e^{-\psi_{m}}d\lambda_{n}$$
have a uniform bound independent of $m$ and $\varepsilon\in(0,\frac{1}{8}B_{0})$.

Therefore the integrals $\int_{ K_0}|u_{v,t_0,\varepsilon}|^{2}e^{-\psi}d\lambda_{n}$
have a uniform bound independent of $\varepsilon\in(0,\frac{1}{8}B_{0})$,
for any given compact set $K_{0}\subset\subset U_0\cap D_v$ containing $o$.

In summary,
we have $|u_{v,t_0,\varepsilon}|^{2}e^{-\psi}$ is integrable near $o$,
and $\bar\partial u_{v,t_0,\varepsilon}=0$ near $o$.
That is to say
$$(u_{v,t_0,\varepsilon},o)\in\mathcal{I}(\psi)_{o}.$$

Let $F_{v,t_0,\varepsilon}:=(1-v'_{t_0,\varepsilon}\circ\psi)F-u_{v,t_0,\varepsilon}$.
By inequality \ref{equ:3.3} and $(u_{v,t_0,\varepsilon},o)\in\mathcal{I}(\psi)_{o}$,
it follows that
$F_{v,t_0,\varepsilon}$ is a holomorphic function on $D_{v}$ satisfying $(F_{v,t_0,\varepsilon}-F,o)\in\mathcal{I}(\psi)_{o}$ and
\begin{equation}
\label{equ:3.5}
\begin{split}
&\int_{ D_v}|F_{v,t_0,\varepsilon}-(1-v'_{t_0,\varepsilon}\circ\psi)F|^{2}d\lambda_{n}
\\\leq&\frac{\mathbf{C}}{e^{A_{t_0}}}\int_{D_v}(v''_{t_0,\varepsilon}\circ\psi)|F|^{2}e^{-\psi}d\lambda_{n}.
\end{split}
\end{equation}

Given $t_0$ and $D_{v}$,
it is clear that $(v''_{t_0,\varepsilon}\circ\psi)|F|^{2}e^{-\psi}$ have a uniform bound on $D_{v}$ independent of $\varepsilon$.

Then the integrals
$\int_{D_v}(v''_{t_0,\varepsilon}\circ\psi)|F|^{2}e^{-\psi}d\lambda_{n}$ have a uniform bound independent of $\varepsilon$,
for any given $t_0$ and $D_v$.

As $|(1-v'_{t_0,\varepsilon}\circ\psi)F|^{2}$ have a uniform bound on $D_{v}$ independent of $\varepsilon$,
it follows that
the integrals $\int_{D_v}|(1-v'_{t_0,\varepsilon}\circ\psi)F|^{2}d\lambda_{n}$ have a uniform bound independent of $\varepsilon$,
for any given $t_0$ and $D_v$.

As
\begin{equation}
\label{}
\begin{split}
&\int_{ D_v}|F_{v,t_0,\varepsilon}|^{2}d\lambda_{n}
\\&\leq
\int_{ D_v}|F_{v,t_0,\varepsilon}-(1-v'_{t_0,\varepsilon}\circ\psi)F|^{2}d\lambda_{n}
+\int_{ D_v}|(1-v'_{t_0,\varepsilon}\circ\psi)F|^{2}d\lambda_{n}
\\&\leq\frac{\mathbf{C}}{e^{A_{t_0}}}\int_{D_v}(v''_{t_0,\varepsilon}\circ\psi)|F|^{2}e^{-\psi}d\lambda_{n}
+\int_{ D_v}|(1-v'_{t_0,\varepsilon}\circ\psi)F|^{2}d\lambda_{n},
\end{split}
\end{equation}
then $\int_{ D_v}|F_{v,t_0,\varepsilon}|^{2}d\lambda_{n}$ have a uniform bound independent of $\varepsilon$.

As $\bar\partial F_{v,t_{0},\varepsilon}=0$ when $\varepsilon\to 0$ and the unit ball of $L^{2}(D_{v})$ is weakly compact,
it follows that the weak limit of some weakly convergent subsequence of
$\{F_{v,t_0,\varepsilon}\}_{\varepsilon}$ gives us a holomorphic function $F_{v,t_0}$ on $\Delta^{n}$.

Then we can also choose a subsequence of the weakly convergent subsequence of $\{F_{v,t_0,\varepsilon}\}_{\varepsilon}$,
such that the chosen sequence is uniformly convergent on any compact subset of $D_v$, denoted by $\{F_{v,t_0,\varepsilon}\}_{\varepsilon}$ without ambiguity.

For any given compact subset $K_{0}$ on $D_v$,
$F_{v,t_0,\varepsilon}$, $(1-v'_{t_0,\varepsilon}\circ\psi)F$ and
$(v''_{t_0,\varepsilon}\circ\psi)|F|^{2}e^{-\psi}$ have
uniform bounds on $K_{0}$ independent of $\varepsilon$.

As the integrals $\int_{ K_0}|u_{v,t_0,\varepsilon}|^{2}e^{-\psi}d\lambda_{n}$
have a uniform bound independent of $\varepsilon\in(0,\frac{1}{8}B_{0})$,
for any given compact set $K_{0}\subset\subset U_0\cap D_v$ containing $o$,
it follows that
$$(F_{v,t_0}-(1-b_{t_0}(\psi))F,o)\in \mathcal{I}(\psi)_{o}.$$

Using the dominated convergence theorem on any compact subset $K$ of
$D_v$ and inequality \ref{equ:3.5}, we obtain
\begin{equation}
\begin{split}
&\int_{K}|F_{v,t_0}-(1-b_{t_0}(\psi))F|^{2}d\lambda_{n}
\\\leq&\frac{\mathbf{C}}{e^{A_{t_0}}}\int_{D_v}(\frac{1}{B_{0}}\mathbb{I}_{\{-t_{0}-B_{0}< t<-t_{0}\}}\circ\psi)|F|^{2}e^{-\psi}d\lambda_{n}.
\end{split}
\end{equation}

It suffices to find $\eta$ and $\phi$ such that
$(\eta+g^{-1})\leq \mathbf{C}e^{-\psi_{m}}e^{-\phi}=\mathbf{C}\tilde{\mu}^{-1}$ on $D_v$.
As $\eta=s(-v_{t_0,\varepsilon}\circ\psi_{m})$ and $\phi=u(-v_{t_0,\varepsilon}\circ\psi_{m})$,
we have $(\eta+g^{-1}) e^{v_{t_0,\varepsilon}\circ\psi_{m}}e^{\phi}=(s+\frac{s'^{2}}{u''s-s''})e^{-t}e^{u}\circ(-v_{t_0,\varepsilon}\circ\psi_{m})$.

Summarizing the above discussion about $s$ and $u$, we are naturally led to a
system of ODEs (see \cite{guan-zhou12,guan-zhou13p,guan-zhou13ap,GZopen-b}):
\begin{equation}
\label{GZ}
\begin{split}
&1).\,\,(s+\frac{s'^{2}}{u''s-s''})e^{u-t}=\mathbf{C}, \\
&2).\,\,s'-su'=1,
\end{split}
\end{equation}
where $t\in[0,+\infty)$, and $\mathbf{C}=1$.

It is not hard to solve the ODE system \ref{GZ} and get $u=-\log(1-e^{-t})$ and
$s=\frac{t}{1-e^{-t}}-1$.
It follows that $s\in C^{\infty}((0,+\infty))$ satisfies $s\geq0$, $\lim_{t\to+\infty}u(t)=0$ and
$u\in C^{\infty}((0,+\infty))$ satisfies $u''s-s''>0$.

As $u=-\log(1-e^{-t})$ is decreasing with respect to $t$,
then
\begin{equation}
\begin{split}
&
\frac{\mathbf{C}}{e^{A_{t_0}}}=\frac{1}{e^{\inf_{t\in(t_{0}+B_{0},t_{0})}u(t)}}
\\&=\sup_{t\in(t_{0}+B_{0},t_{0})}\frac{1}{e^{u(t)}}
=\sup_{t\in(t_{0}+B_{0},t_{0})}(1-e^{-t})=1-e^{-(t_{0}+B_{0})},
\end{split}
\end{equation}
therefore we are done.
Thus we prove Lemma \ref{p:GZ_JM_sharp}.

\subsection{Proof of Proposition \ref{p:GZ_JM201312}}
$\\$

Let
$$\varphi:=(1+\delta)\max\{\psi,\log|F|^{2}\}.$$

It suffices to prove the case
that
$$\Psi:=\min\{\psi-\log|F|^{2},0\}-a,$$
where $a>0$ is arbitrarily given.

For the sake of completeness, we recall our proof in \cite{guan-zhou13p} (see also \cite{guan-zhou13ap}) with slightly modifications.

As $D_{v}\subset\subset\Delta^{n}\subset \mathbb{C}^{n}$,
then there exist negative smooth plurisubharmonic functions $\{\varphi_{m}\}_{m=1,2,\cdots}$ and
smooth functions $\{\Psi_{m}\}_{m=1,2,\cdots}$
on a neighborhood of $\overline{D}_{v}$,
such that

1). $\{\varphi_{m}+\Psi_{m}\}_{m=1,2,\cdots}$ and $\{\varphi_{m}+(1+\delta)\Psi_{m}\}_{m=1,2,\cdots}$ are negative smooth plurisubharmonic functions;

2). the sequence $\{\varphi_{m}\}_{m=1,2,\cdots}$ is decreasingly convergent to $\varphi$;

3). the sequence $\{\varphi_{m}+\Psi_{m}\}_{m=1,2,\cdots}$ is decreasingly convergent to $\varphi+\Psi$;
$\\$on a smaller neighborhood of $\overline{D}_{v}$,
when $m\to+\infty$.

Let $\eta=s(-v_{t_{0},\varepsilon}\circ\Psi_{m})$ and $\phi=u(-v_{t_{0},\varepsilon}\circ\Psi_{m})$,
where $s\in C^{\infty}((0,+\infty))$ satisfies $s\geq\frac{1}{\delta}$, and
$u\in C^{\infty}((0,+\infty))$, such that $u''s-s''>0$, and $s'-u's=1$.

Let $\Phi:=\varphi_{m}+\Psi_{m}+\phi$.

Now let $\alpha=\sum^{n}_{j=1}\alpha_{j}d\bar z^{j}\in Dom_{D_v}
(\bar\partial^*)\cap Ker(\bar\partial)\cap C^\infty_{(0,1)}(\overline {D_v})$.
By Cauchy-Schwarz inequality, it follows that
\begin{equation}
\label{equ:20131130aJM}
\begin{split}
2\mathrm{Re}(\bar\partial^*_\Phi\alpha,\alpha\llcorner(\bar\partial\eta)^\sharp )_{\Omega,\Phi}
\geq
&-\int_{D_v}g^{-1}|\bar{\partial}^{*}_{\Phi}\alpha|^{2}e^{-\Phi}d\lambda_{n}
\\&+
\sum_{j,k=1}^{n}\int_{D_v}
(-g(\partial_{j} \eta)\bar\partial_{k} \eta )\alpha_{\bar{j} }\overline{{\alpha}_{\bar {k}}}e^{-\Phi}d\lambda_{n}.
\end{split}
\end{equation}

Using the twisted Bochner-Kodaira identity (see Lemma 3.2. in \cite{guan-zhou-zhu10}) and inequality \ref{equ:20131130aJM},
since $s\geq0$ and $\varphi_{m}$ is a plurisubharmonic function on $\overline{D}_{v}$,
we get

\begin{equation}
\label{equ:4.1JM}
\begin{split}
&\int_{D_v}(\eta+g^{-1})|\bar{\partial}^{*}_{\Phi}\alpha|^{2}e^{-\Phi}d\lambda_{n}
\\\geq&\sum_{j,k=1}^{n}\int_{D_v}
(-\partial_{j}\bar{\partial}_{k}\eta+\eta\partial_{j}\bar{\partial}_{k}\phi+\eta\partial_{j}\bar{\partial}_{k}(\Psi_{m}+\varphi_{m})-g(\partial_{j} \eta)\bar\partial_{k} \eta )\alpha_{\bar{j} }\overline{{\alpha}_{\bar {k}}}e^{-\Phi}d\lambda_{n},
\end{split}
\end{equation}
where $g$ is a positive continuous function on $D_{v}$.

We need some calculations to determine $g$.

We have
\begin{equation}
\label{equ:20131204bJM}
\begin{split}
&\sum_{1\leq j,k\leq n}(-\partial_{j}\bar{\partial}_{k}\eta+\eta\partial_{j}\bar{\partial}_{k}\phi-g(\partial_{j} \eta)
\bar\partial_{k} \eta )\alpha_{\bar{j} }\overline{{\alpha}_{\bar{k}}}
\\=&(s'-su')\sum_{1\leq j,k\leq n}\partial_{j}\bar{\partial}_{k}(v_{t_0,\varepsilon}\circ \Psi_{m})\alpha_{\bar{j} }\overline{{\alpha}_{\bar{k}}}
\\+&((u''s-s'')-gs'^{2})\sum_{1\leq j,k\leq n}\partial_{j}
(-v_{t_0,\varepsilon}\circ \Psi_{m})\bar{\partial}_{k}(-v_{t_0,\varepsilon}\circ\Psi_{m})\alpha_{\bar{j} }\overline{{\alpha}_{\bar{k}}}
\\=&(s'-su')\sum_{1\leq j,k\leq n}((v'_{t_0,\varepsilon}\circ\Psi_{m})\partial_{j}\bar{\partial}_{k}\Psi_{m}+(v''_{t_0,\varepsilon}\circ \Psi_{m})\partial_{j}(\Psi_{m})\bar{\partial}_{k}(\Psi_{m}))\alpha_{\bar{j} }\overline{{\alpha}_{\bar{k}}}
\\+&((u''s-s'')-gs'^{2})\sum_{1\leq j,k\leq n}\partial_{j}
(-v_{t_0,\varepsilon}\circ \Psi_{m})\bar{\partial}_{k}(-v_{t_0,\varepsilon}\circ \Psi_{m})\alpha_{\bar{j} }\overline{{\alpha}_{\bar{k}}}.
\end{split}
\end{equation}
We omit composite item $(-v_{t_0,\varepsilon}\circ \Psi_{m})$ after $s'-su'$ and $(u''s-s'')-gs'^{2}$ in the above equalities.

Since $\varphi_{m}+\Psi_{m}$ and $\varphi_{m}+(1+\delta)\Psi_{m}$ are plurisubharmonic on $\overline{D}_{v}$ and
$0\leq v'_{t_{0},\varepsilon}\circ\Psi_{m}\leq1$,
we have
\begin{equation}
(1-v'_{t_0,\varepsilon}\circ\Psi_{m})\sqrt{-1}\partial\bar\partial(\varphi_{m}+\Psi_{m})+
(v'_{t_0,\varepsilon}\circ\Psi_{m})\sqrt{-1}\partial\bar\partial(\varphi_{m}+(1+\delta)\Psi_{m})\geq 0,
\end{equation}
on $\overline{D}_{v}$, which means that
\begin{equation}
\label{equ:20131204aJM}
\frac{1}{\delta}\sqrt{-1}\partial\bar\partial(\varphi_{m}+\Psi_{m})+(v'_{t_0,\varepsilon}
\circ\Psi_{m})\sqrt{-1}\partial\bar{\partial}\Psi_{m}\geq 0,
\end{equation}
on $\overline{D}_{v}$.

Let $g=\frac{u''s-s''}{s'^{2}}\circ(-v_{t_0,\varepsilon}\circ \Psi_{m})$.
It follows that $\eta+g^{-1}=(s+\frac{s'^{2}}{u''s-s''})\circ(-v_{t_0,\varepsilon}\circ \Psi_{m})$.

Because of $v'_{t_0,\varepsilon}\geq 0$  and $s'-su'=1$, using inequalities \ref{equ:4.1JM} \ref{equ:20131204aJM} and \ref{equ:20131204bJM},
we have
\begin{equation}
\label{equ:3.1JM}
\begin{split}
\int_{D_v}(\eta+g^{-1})|\bar{\partial}^{*}_{\Phi}\alpha|^{2}e^{-\Phi}d\lambda_{n}
\geq\int_{D_v}(v''_{t_0,\varepsilon}\circ{\Psi_{m}})
\big|\alpha\llcorner(\bar \partial \Psi_{m})^\sharp\big|^2e^{-\Phi}d\lambda_{n}.
\end{split}
\end{equation}

Let $\lambda=\bar{\partial}[(1-v'_{t_0,\varepsilon}(\Psi_{m})){F^{1+\delta}}]$.
By the definition of contraction, Cauchy-Schwarz inequality and inequality \ref{equ:3.1JM},
it follows that
\begin{equation}
\begin{split}
&|(\lambda,\alpha)_{D_v,\Phi}|^{2}=|((v''_{t_0,\varepsilon}\circ{\Psi_{m}})F^{1+\delta},\alpha\llcorner(\bar\partial\Psi_{m})^\sharp\big)_{D_v,\Phi}|^{2}
\\\leq&
(\int_{D_v}
(v''_{t_0,\varepsilon}\circ{\Psi_{m}})|F^{1+\delta}|^2e^{-\Phi}d\lambda_{n})
(\int_{ D_v}(\eta+g^{-1})|\bar{\partial}^{*}_{\Phi}\alpha|^{2}e^{-\Phi}d\lambda_{n}).
\end{split}
\end{equation}

Let $\mu:=(\eta+g^{-1})^{-1}$. Using a lemma of Berndtsson's (see \cite{berndtsson} or Lemma 3.7 in \cite{guan-zhou-zhu10}),
we have locally $L^{1}$ function $u_{v,t_0,m,\varepsilon}$ on $D_{v}$ such that $\bar{\partial}u_{v,t_0,m,\varepsilon}=\lambda$,
and
\begin{equation}
 \label{equ:3.2JM}
 \begin{split}
 &\int_{ D_v}|u_{v,t_0,m,\varepsilon}|^{2}(\eta+g^{-1})^{-1} e^{-\Phi}d\lambda_{n}
  \leq\int_{D_v}(v''_{t_0,\varepsilon}\circ{\Psi_{m}})|F^{1+\delta}|^2e^{-\Phi}d\lambda_{n}.
  \end{split}
\end{equation}

Let $\mu_{1}=e^{v_{t_0,\varepsilon}\circ\Psi_{m}}$, $\tilde{\mu}=\mu_{1}e^{\phi}$.
Assume that we can choose $\eta$ and $\phi$ such that $\tilde{\mu}\leq \mathbf{C}(\eta+g^{-1})^{-1}=\mu$, where $\mathbf{C}=1$.

Note that $v_{t_0,\varepsilon}(\Psi_{m})\geq\Psi_{m}$.
Then it follows that
\begin{equation}
\label{equ:3.8JM}
\begin{split}
\int_{ D_v}|u_{v,t_0,m,\varepsilon}|^{2}e^{-\varphi_{m}} d\lambda_{n}
&\leq\int_{ D_v}|u_{v,t_0,m,\varepsilon}|^{2}\tilde{\mu}e^{-\Phi}d\lambda_{n}.
\end{split}
\end{equation}

Using inequalities \ref{equ:3.2JM} and \ref{equ:3.8JM},
we obtain that
$$\int_{ D_v}|u_{v,t_0,m,\varepsilon}|^{2}e^{-\varphi_{m}} d\lambda_{n}\leq\mathbf{C}\int_{D_v}
(v''_{t_0,\varepsilon}\circ\Psi_{m})|F^{1+\delta}|^2e^{-\Phi}d\lambda_{n},$$
under the assumption $\tilde{\mu}\leq\mathbf{C} (\eta+g^{-1})^{-1}$.

As $-v_{t_0,\varepsilon}\circ{\Psi_{m}}(\overline{D_{v}})\subset\subset(-\infty,t_{0}+1)$,
then it is clear that
\begin{equation}
\label{equ:20131130bJM}
-v_{t_0,\varepsilon}\circ{\Psi_{m}}(\overline{D_{v}})\subset(-\infty,K_{t_{0}})
\end{equation}
where $K_{t_{0}}$ is independent of $m$ and $\varepsilon\in (0,\frac{1}{8}B_{0})$.

As $u$ is positive and smooth on $(-\infty,+\infty)$,
it follows that $\phi$ is uniformly bounded on $\overline{D}_{v}$ independent of $m$.

As $Supp(v''_{t_{0},\varepsilon})\subset\subset(-t_{0}-B_{0},-t_{0})$
and
$$|F^{1+\delta}|^2e^{-\varphi_{m}}\leq|F^{1+\delta}|^2e^{-\varphi}\leq e^{-(1+\delta)\max\{\psi-\log|F|^{2},0\}}\leq 1,$$
then it is clear that $(v''_{t_0,\varepsilon}\circ{\Psi_{m}})|F^{1+\delta}|^2e^{-\varphi_{m}-\Psi_{m}}$ are uniformly bounded on $\overline{D}_{v}$ independent of $m$.

Therefore the integrals
$\int_{D_v}(v''_{t_0,\varepsilon}\circ{\Psi_{m}})|F^{1+\delta}|^2e^{-\Phi}d\lambda_{n}$ are uniformly bounded independent of $m$,
for any given $v$, $t_0$, $\varepsilon$.

By weakly compactness of the unit ball of $L^{2}_{\varphi}(D_{v})$ and dominated convergence theorem, when $m\to+\infty$,
it follows that the weak limit of some weakly convergent subsequence of
$\{u_{v,t_{0},m,\varepsilon}\}_{m}$ gives function $u_{v,t_0,\varepsilon}$ on $D_v$ satisfying
\begin{equation}
\label{equ:3.3JM}
\int_{ D_v}|u_{v,t_0,\varepsilon}|^{2}e^{-\varphi}d\lambda_{n}\leq\frac{\mathbf{C}}{e^{A_{t_0}}}\int_{D_v}
(v''_{t_0,\varepsilon}\circ{\Psi})|F^{1+\delta}|^2e^{-\varphi-\Psi}d\lambda_{n},
\end{equation}
where $A_{t_0}:=\inf_{t\geq t_0}\{u(t)\}$.

Let $F_{v,t_0,m,\varepsilon}:=(1-v'_{t_0,\varepsilon}(\Psi_{m}))F^{1+\delta}-u_{v,t_0,m,\varepsilon}$,
which is a holomorphic function on $D_{v}$.

As $|F^{1+\delta}|^{2}e^{-\varphi_{m}}\leq|F^{1+\delta}|^{2}e^{-\varphi} \leq 1$,
then the integrals
$$\int_{ D_v}|(1-v'_{t_0,\varepsilon}(\Psi_{m}))F^{1+\delta}|^{2}e^{-\varphi_{m}} d\lambda_{n}$$
have a uniform bound independent of $m$.
Recall that the integrals
$$\int_{ D_v}|u_{v,t_{0},m,\varepsilon}|^{2}e^{-\varphi_{m}} d\lambda_{n}$$
have a uniform bound independent of $m$,
then the integrals
$$\int_{ D_v}|F_{v,t_0,m,\varepsilon}|^{2}e^{-\varphi_{m}} d\lambda_{n}$$
have a uniform bound independent of $m$.

Therefore we can choose a subsequence of $\{F_{v,t_0,m,\varepsilon}\}_{m=1,2,\cdots}$ from the chosen weakly convergent subsequence of $$(1-v'_{t_0,\varepsilon}(\Psi_{m}))F^{1+\delta}-u_{v,t_0,m,\varepsilon},$$
such that the subsequence is uniformly convergent on any compact subset of $D_{v}$,
and we still denote the subsequence by $\{F_{v,t_0,m,\varepsilon}\}_{m=1,2,\cdots}$ without ambiguity.

Denote by
$$F_{v,t_0,\varepsilon}:=\lim_{m\to\infty}F_{v,t_0,m,\varepsilon}.$$

By the above arguments,
it follows that the right hand side of inequality \ref{equ:3.2JM} are uniformly bounded independent of $m$ and $\varepsilon\in(0,\frac{1}{8}B_{0})$.

By inequality \ref{equ:3.2JM},
it follows that the integrals
$$\int_{ D_v}|u_{v,t_0,m,\varepsilon}|^{2}(\eta+g^{-1})^{-1} e^{-\phi-\varphi_{m}-\Psi_{m}}d\lambda_{n}$$
are uniformly bounded independent of $m$ and $\varepsilon\in(0,\frac{1}{8}B_{0})$.

Using inequality \ref{equ:20131130bJM},
we obtain that
$$(\eta+g^{-1})^{-1}=(s(-v_{t_{0},\varepsilon}\circ\Psi_{m})+\frac{s'^{2}}{u''s-s''}\circ(-v_{t_0,\varepsilon}\circ\Psi_{m}))^{-1}$$
and
$e^{-\phi}=e^{-u(-v_{t_{0},\varepsilon}\circ\Psi_{m})}$ have positive uniform lower bounds independent of $m$ and $\varepsilon\in(0,\frac{1}{8}B_{0})$.

Then the integrals
$$\int_{ K_0}|u_{v,t_0,m,\varepsilon}|^{2}e^{-\varphi_{m}-\Psi_{m}}d\lambda_{n}$$
have a uniform upper bound independent of $m$ and $\varepsilon\in(0,\frac{1}{8}B_{0})$
for any given compact set $K_{0}\subset\subset D_v$.

As $$Supp(v'_{t_0,\varepsilon}(\Psi_{m}))\subset \{\Psi_{m}>-t_{0}-1\},$$
it follows that
$$|v'_{t_0,\varepsilon}(\Psi_{m})|^{2}e^{-\Psi_{m}}\leq e^{t_{0}+1}.$$

Furthermore, as
$$|F^{1+\delta}|^{2}e^{-\varphi_{m}}\leq|F^{1+\delta}|^{2}e^{-\varphi}=e^{-(1+\delta)\max\{\psi-\log|F|^{2},0\}}\leq1$$
then the integrals
$$\int_{ K_0}|v'_{t_0,\varepsilon}(\Psi_{m})F^{1+\delta}|^{2}e^{-\varphi_{m}-\Psi_{m}}d\lambda_{n}$$
have a uniform upper bound independent of $m$ and $\varepsilon\in(0,\frac{1}{8}B_{0})$
for any given compact set $K_{0}\subset\subset D_v$.
Therefore the integrals
\begin{equation}
\label{equ:20131211JM}\int_{ K_0}|F_{v,t_0,m,\varepsilon}-F^{1+\delta}|^{2}e^{-\varphi_{m}-\Psi_{m}}d\lambda_{n}
\end{equation}
have a uniform upper bound independent of $m$ and $\varepsilon\in(0,\frac{1}{8}B_{0})$
for any given compact set $K_{0}\subset\subset D_v$.

As
$\varphi_{m'}+\Psi_{m'}\leq\varphi_{m}+\Psi_{m}$ where $m'\geq m$,
it follows that
$$|F_{v,t_0,m',\varepsilon}-F^{1+\delta}|^{2}e^{-(\varphi_{m}+\Psi_{m})}\leq|F_{v,t_0,m',\varepsilon}-F^{1+\delta}|^{2}e^{-(\varphi_{m'}+\Psi_{m'})}.$$

By inequality \ref{equ:20131211JM}, it follows that for any given compact set $K_{0}\subset\subset  D_v$,
the integrals
$$\int_{ K_0}|F_{v,t_0,m',\varepsilon}-F^{1+\delta}|^{2}e^{-(\varphi_{m}+\Psi_{m})}d\lambda_{n}$$
have a uniform bound independent of $m$ and $m'$.

Therefore
for any given compact set $K_{0}\subset\subset\cap D_v$,
the integrals
$$\int_{ K_0}|F_{v,t_0,\varepsilon}-F^{1+\delta}|^{2}e^{-(\varphi_{m}+\Psi_{m})}d\lambda_{n}$$
have a uniform bound independent of $m$ and $\varepsilon\in(0,\frac{1}{8}B_{0})$.

It is clear that the integrals
\begin{equation}
\label{equ:20131211aJM}\int_{ K_0}|F_{v,t_0,\varepsilon}-F^{1+\delta}|^{2}e^{-\varphi-\Psi}d\lambda_{n}
\end{equation}
have a uniform upper bound independent of $\varepsilon\in(0,\frac{1}{8}B_{0})$
for any given compact set $K_{0}\subset\subset D_v$.

In summary,
we have $|F_{v,t_0,\varepsilon}-F^{1+\delta}|^{2}e^{-\varphi-\Psi}$ is integrable near $o$.
That is to say
$$(F_{v,t_0,\varepsilon}-F^{1+\delta},o)\in\mathcal{I}(\varphi+\Psi)_{o}.$$

By inequality \ref{equ:3.3JM},
it follows that
$F_{v,t_0,\varepsilon}$ is a holomorphic function on $D_{v}$ satisfying $(F_{v,t_0,\varepsilon}-F^{1+\delta},o)\in\mathcal{I}(\varphi)_{o}$ and
\begin{equation}
\label{equ:3.5JM}
\begin{split}
&\int_{ D_v}|F_{v,t_0,\varepsilon}-(1-v'_{t_0,\varepsilon}\circ\Psi)F^{1+\delta}|^{2}e^{-\varphi}d\lambda_{n}
\\\leq&\frac{\mathbf{C}}{e^{A_{t_0}}}\int_{D_v}(v''_{t_0,\varepsilon}\circ\Psi)|F^{1+\delta}|^{2}e^{-\varphi-\Psi}d\lambda_{n}.
\end{split}
\end{equation}

Given $t_0$ and $D_{v}$,
it is clear that $(v''_{t_0,\varepsilon}\circ\Psi)|F^{1+\delta}|^{2}e^{-\varphi-\Psi}$ have a uniform bound on $D_{v}$ independent of $\varepsilon$.
Then the integrals
$\int_{D_v}(v''_{t_0,\varepsilon}\circ\Psi)|F^{1+\delta}|^{2}e^{-\varphi-\Psi}d\lambda_{n}$ have a uniform bound independent of $\varepsilon$,
for any given $t_0$ and $D_v$.

As $|(v'_{t_0,\varepsilon}\circ\Psi)F^{1+\delta}|^{2}e^{-\varphi}$ have a uniform bound on $D_{v}$ independent of $\varepsilon$,
it follows that the integrals $\int_{D_v}|(1-v'_{t_0,\varepsilon}\circ\Psi)F^{1+\delta}|^{2}e^{-\varphi}d\lambda_{n}$ have a uniform bound independent of $\varepsilon$,
for any given $t_0$ and $D_v$.

As
\begin{equation}
\label{}
\begin{split}
&\int_{ D_v}|F_{v,t_0,\varepsilon}|^{2}e^{-\varphi}d\lambda_{n}
\\&\leq
\int_{ D_v}|F_{v,t_0,\varepsilon}-(1-v'_{t_0,\varepsilon}\circ\Psi)F^{1+\delta}|^{2}e^{-\varphi}d\lambda_{n}
+\int_{ D_v}|(1-v'_{t_0,\varepsilon}\circ\Psi)F^{1+\delta}|^{2}e^{-\varphi}d\lambda_{n}
\\&\leq\frac{\mathbf{C}}{e^{A_{t_0}}}\int_{D_v}(v''_{t_0,\varepsilon}\circ\Psi)|F^{1+\delta}|^{2}e^{-\varphi-\Psi}d\lambda_{n}
+\int_{ D_v}|(1-v'_{t_0,\varepsilon}\circ\Psi)F^{1+\delta}|^{2}e^{-\varphi}d\lambda_{n},
\end{split}
\end{equation}
then the integrals $\int_{ D_v}|F_{v,t_0,\varepsilon}|^{2}e^{-\varphi}d\lambda_{n}$ have a uniform bound independent of $\varepsilon$.

As $\bar\partial F_{v,t_{0},\varepsilon}=0$ when $\varepsilon\to 0$ and the unit ball of $L^{2}_{\varphi}(D_{v})$ is weakly compact,
it follows that the weak limit of some weakly convergent subsequence of
$\{F_{v,t_0,\varepsilon}\}_{\varepsilon}$ gives us a holomorphic function $F_{v,t_0}$ on $\Delta^{n}$.

Then we can also choose a subsequence of the weakly convergent subsequence of $\{F_{v,t_0,\varepsilon}\}_{\varepsilon}$,
such that the chosen sequence is uniformly convergent on any compact subset of $D_v$, denoted by $\{F_{v,t_0,\varepsilon}\}_{\varepsilon}$ without ambiguity.

For any given compact subset $K_{0}$ on $D_v$, $F_{v,t_0,\varepsilon}$,
$|(1-v'_{t_0,\varepsilon}\circ\Psi)F^{1+\delta}|^{2}e^{-\varphi}$ and $(v''_{t_0,\varepsilon}\circ\Psi)|F^{1+\delta}|^{2}e^{-\varphi-\Psi}$
have uniform bounds on $K_{0}$ independent of $\varepsilon$.

By inequality \ref{equ:20131211aJM},
it follows that
$$(F_{v,t_0}-F^{1+\delta},o)\in \mathcal{I}(\varphi+\Psi)_{o}.$$

Using the dominated convergence theorem on any compact subset $K$ of $D_v$,
we obtain
\begin{equation}
\begin{split}
&\int_{K}|F_{v,t_0}-(1-b_{t_0}(\Psi))F^{1+\delta}|^{2}e^{-\varphi}d\lambda_{n}
\\\leq&(1+\frac{1}{\delta})\int_{D_v}\frac{1}{B_{0}}(\mathbb{I}_{\{-t_{0}-B_{0}< t<-t_{0}\}}\circ\Psi)|F^{1+\delta}|^{2}e^{-\varphi-\Psi}d\lambda_{n}.
\end{split}
\end{equation}

Proposition \ref{p:GZ_JM201312} has thus been proved.

It suffices to find $\eta$ and $\phi$ such that
$(\eta+g^{-1})\leq \mathbf{C}e^{-\Psi_{m}}e^{-\phi}=\mathbf{C}\tilde{\mu}^{-1}$ on $D_v$.
As $\eta=s(-v_{t_0,\varepsilon}\circ\Psi_{m})$ and $\phi=u(-v_{t_0,\varepsilon}\circ\Psi_{m})$,
we have $(\eta+g^{-1}) e^{v_{t_0,\varepsilon}\circ\Psi_{m}}e^{\phi}=(s+\frac{s'^{2}}{u''s-s''})e^{-t}e^{u}\circ(-v_{t_0,\varepsilon}\circ\Psi_{m})$.

Summarizing the above discussion about $s$ and $u$, we are naturally led to a
system of ODEs:
\begin{equation}
\label{GZJM}
\begin{split}
&1).\,\,(s+\frac{s'^{2}}{u''s-s''})e^{u-t}=\mathbf{C}, \\
&2).\,\,s'-su'=1,
\end{split}
\end{equation}
where $t\in[0,+\infty)$, and $\mathbf{C}=1$.

It is not hard to solve the ODE system \ref{GZJM} (using the same arguments in \cite{guan-zhou13p} \cite{guan-zhou13ap} and noting the boundary condition)
and get $u=-\log(1+\frac{1}{\delta}-e^{-t})$ and
$s=\frac{(1+\frac{1}{\delta})t+\frac{1}{\delta}(1+\frac{1}{\delta})}{1+\frac{1}{\delta}-e^{-t}}-1$.
It follows that $s\in C^{\infty}((0,+\infty))$ satisfies $s\geq\frac{1}{\delta}$, $u'\leq 0$ and
$u\in C^{\infty}((0,+\infty))$ satisfies $u''s-s''>0$.

As $u=-\ln(1+\frac{1}{\delta}-e^{-t})$ is decreasing with respect to $t$,
then
$$\frac{\mathbf{C}}{e^{A_{t_0}}}=\frac{1}{\exp\inf_{t\geq t_{0}}u(t)}=\sup_{t\geq t_{0}}\frac{1}{e^{u(t)}}=\sup_{t\geq t_{0}}(1+\frac{1}{\delta}-e^{-t})=1+\frac{1}{\delta},$$
for any $t_{0}\geq0$,
therefore we are done.

\section{Discussion of inequality \ref{equ:sharp20140126}}\label{sec:2014spring}

We would like to give a proof of and a Remark on inequality \ref{equ:sharp20140126}.

\subsection{Proof of inequality \ref{equ:sharp20140126}}\label{sec:2014springb}
$\\$

It suffices to prove
\begin{equation}
\label{equ:sharp20140126b}
\frac{t}{6(t-1)}<(\frac{1}{(t-1)(2t-1)})^{\frac{1}{t}}.
\end{equation}

We consider the function
$$P(t)=\frac{1}{t}\log\frac{1}{(t-1)(2t-1)}+\log\frac{t-1}{t}.$$

Replacing $t$ by $\frac{1}{x}$,
we obtain
\begin{equation}
\begin{split}
Q(x)=P(\frac{1}{x})
&=\frac{1}{\frac{1}{x}}\log\frac{1}{(\frac{1}{x}-1)(2\frac{1}{x}-1)}+\log\frac{\frac{1}{x}-1}{\frac{1}{x}}
\\&=\log((\frac{x^{2}}{(1-x)(2-x)})^{x}(1-x))
\\&=2x\log x+(1-x)\log(1-x)-x\log(2-x)
\end{split}
\end{equation}

We need to prove that inequality $e^{Q(x)}>\frac{1}{6}$ holds for any $x\in(0,1)$.

One can obtain the derivative $Q'(x)$ of $Q(x)$ as follows

\begin{equation}
\begin{split}
Q'(x)&=(2x\log x+(1-x)\log(1-x)-x\log(2-x))'
\\&=(2\log x+2)+(-1-\log(1-x))+(\frac{x}{2-x}-\log(2-x))
\\&=2\log x-\log(1-x)+\frac{2}{2-x}-\log(2-x)
\end{split}
\end{equation}

One can obtain the derivative $Q''(x)$ of $Q(x)$ as follows

\begin{equation}
\begin{split}
Q''(x)&=(2\log x-\log(1-x)+\frac{2}{2-x}-\log(2-x))'
\\&=\frac{2}{x}+\frac{1}{1-x}+\frac{2}{(2-x)^{2}}+\frac{1}{2-x}>0
\end{split}
\end{equation}

Note that
$$\lim_{x\to0}Q(x)=0,$$
$$\lim_{x\to1}Q(x)=0,$$
and $Q'(x)(\frac{1}{2})>0$,
then the preimage of the minimal value must exist in $(0,\frac{1}{2})$.

Now we consider the minimal of
the following three functions:

$Q_{1}(x):=2x\log x$, $Q_{2}(x):=(1-x)\log(1-x)$, $Q_{3}(x):=-x\log(2-x)$.

It is clear that
$$Q(x)=Q_{1}(x)+Q_{2}(x)+Q_{3}(x),$$
and
$$\min_{x\in(0,\frac{1}{2})}Q(x)\geq \min_{x\in(0,\frac{1}{2})}Q_{1}(x)+\min_{x\in(0,\frac{1}{2})}Q_{2}(x)+\min_{x\in(0,\frac{1}{2})}Q_{3}(x).$$

It is known that $\min_{x\in(0,\frac{1}{2})}Q_{1}(x)=-\frac{2}{e}$,
$\min_{x\in(0,\frac{1}{2})}Q_{2}(x)>-\log\sqrt{2}$,
and $\min_{x\in(0,\frac{1}{2})}Q_{3}(x)=Q_{3}(\frac{1}{2})>-\log\sqrt{\frac{3}{2}}$.
Then we obtain that
$$\min_{x\in(0,1)}Q(x)\geq\min_{x\in(0,\frac{1}{2})}Q(x)> -\frac{2}{e}-\log\sqrt{2}-\log\sqrt{\frac{3}{2}},$$
that is to say
\begin{equation}
\label{equ:sharp20140126c}
\min_{x\in(0,1)}e^{Q(x)}>\frac{1}{\sqrt{3}e^{\frac{2}{e}}}>\frac{1}{6}.
\end{equation}

Thus we obtain inequality \ref{equ:sharp20140126}.

\subsection{A Remark on inequality \ref{equ:sharp20140126}}\label{sec:2014spring}
$\\$

In this subsection,
we give a remark about the accuracy our effectiveness of Corollary \ref{c:effect}.

When $D=\Delta$,
and $z_{0}=0$, $\varphi=\frac{1}{p}\log|z|^{2}$,
it is clear that $\parallel 1\parallel^{2}_{\varphi}K(z_{0})=\frac{1}{1-\frac{1}{p}}$.

Using Corollary \ref{c:effect} and inequality \ref{equ:sharp20140126c},
we obtain
$$\parallel 1\parallel^{2}_{\varphi}K(z_{0})> \frac{1}{\sqrt{3}e^{\frac{2}{e}}}\frac{1}{1-\frac{1}{p}}$$
for any $D$ and any $z_{0}\in D$.

One can obtain that
$\min_{x\in(0,1)}e^{Q(x)}>0.2876,$
which gives a more precise form of inequality \ref{equ:sharp20140126}:
\begin{equation}
\frac{1}{400(t-1)}<0.2876\frac{t}{(t-1)}<(\frac{1}{(t-1)(2t-1)})^{\frac{1}{t}}.
\end{equation}
Then we obtain
$$\parallel 1\parallel^{2}_{\varphi}K(z_{0})>0.2876\frac{1}{1-\frac{1}{p}}$$
for any $D$ and any $z_{0}\in D$.

\vspace{.1in} {\em Acknowledgements}. The authors would like to
thank Prof. Bo Berndtsson, Prof. Jean-Pierre Demailly, Prof. Nessim Sibony,
and Prof. Yum-Tong Siu for giving series of talks at CAS and
explaining us their related works.

\bibliographystyle{references}
\bibliography{xbib}

\end{document}